\documentclass[12pt,a4paper]{article}

\usepackage{graphicx}
\usepackage{amsmath,amsfonts}
\usepackage{mathrsfs}
\usepackage{bbold}

\let\bb\mathbb       
\newtheorem{theorem}{Theorem}

\newtheorem{proposition}{Proposition}

\let\scr\mathscr

\def \Liminf{\mathop{\underline{\lim}}\limits}
\def \Limsup{\mathop{\overline{\lim}}\limits}
\def\zs#1{_{\lower 3pt \hbox{$\scriptstyle#1$}}}
\def\Pb{\mathbf{P}}

\def\Ex{\mathbf{E}}
\DeclareMathOperator{\argmax}{argmax}

\def\DD{{\bb D}}
\def\KK{{\bb K}}

\def\WW{{\bb W}}
\def\VV{{\bb V}}\def\BB{{\bb B}}
\def\RR{{\scr R}}\def\DD{{\bb D}}
\def\UU{{\bb U}}    
\def\UU{\mathbf{U}}

\def\sgn{{\rm sgn}} 
\def\1{\mbox{1\hspace{-.25em}I}}  

\let\scr\mathscr

\begin{document}
\title{On Identification of the Threshold  Diffusion Processes}
\author{Yury A. Kutoyants\\
{\small Laboratoire de Statistique et Processus, Universit\'e du Maine}\\
 {\small 72085 Le Mans, C\'edex 9, France  }}

\date{}

\maketitle
\begin{abstract}
We consider the problems of parameter estimation for several models of
threshold ergodic diffusion processes in the asymptotics of large samples.
These models are the direct continuous time analogues of the well-known in
time series analysis threshold autoregressive (TAR) models. In such models the
trend is switching when the observed process atteints some (unknown) values
and the problem is to estimate it or to test some hypotheses concerning these
values. The related statistical problems correspond to the singular estimation
or testing, for example, the rate of convergence of estimators is $T$ and not
$\sqrt{T}$ as in regular estimation problems. We study the asymptotic behavior
of the maximum likelihood and bayesian  estimators and discuss the possibility
of the construction of the goodness of fit test for such models of observation.
\end{abstract}


\section{Introduction}

The simplest example of the threshold  model is the following threshold 
autoregressive (TAR) time series:
\begin{equation}
\label{1}
X_{j+1}=\varrho_1 \,X_j\,\1_{\left\{X_j<\vartheta
\right\}}+\varrho_2 \,X_j\,\1_{\left\{X_j\geq \vartheta
\right\}}+\varepsilon _{j+1}, \qquad j=0,\ldots, n-1,
\end{equation}
where $\varepsilon _j$ are i.i.d. ${\cal N}\left(0,s^2 \right)$, $\varrho
_1\not =\varrho _2$ and $\left|\varrho_i\right| <1$. Therefore we have two
different autoregressive processes depending on the region of observations
$\left\{x:\;x<\vartheta \right\}$ or $\left\{x:\;x\geq \vartheta
\right\}$. This time series has ergodic properties with invariant density
close to a weighted sum of two Gaussian densities.  If we suppose that
$s^2,\varrho _1,\varrho _2 $ are known and $\vartheta \in \Theta =\left(\alpha
,\beta \right)$ is unknown parameter, then we obtain the first problem of
threshold $\vartheta $ estimation.  It is easy to see that the likelihood
ratio is a piece wise constant (discontinuous) function of $\vartheta $, the
Fisher information is equal infinity. As usual in singular estimation
problems, the rate of convergence of maximum likelihood $\hat\vartheta _n$ or
Bayesian $\tilde\vartheta _n$ estimators is $n$ and not $\sqrt{n}$ i.e.; the
quantities $n\left(\hat\vartheta _n-\vartheta \right) $ and
$n\left(\tilde\vartheta _n-\vartheta \right) $ have non degenerate limits.

 There are many different threshold regression  models of such type  
extensively developed in econometrics and, of course, the identification of
these models attracts attention of statisticians (see, e.g. the works by
Quandt (1958), Tong (1990) \cite{To90}, Chan (1993) \cite{Ch93}, Hansen (2000)
\cite{Han00}, Fan and Yao (2003) \cite{FY03}, Koul {\it et al.} \cite{KLS},
Chan and Kutoyants \cite{CK08} and the references therein). Note that
continuous time models actually  find a wide range of applications in
econometrical problems and occupy a central place in financial mathematics
(see, e.g., the work by Shreve \cite{S}). 

Our goal is to study several models of continuous time analogues
(diffusion processes) of such threshold type time series and to describe the
properties of estimators of the thresholds for these models. Note that the
general theory of parameter estimation (in regular case) for ergodic diffusion
processes is actually well developped (see, e.g. \cite{Kut04}, \cite{NY} and
references therein) but the problems of threshold estimation are of singular
type and need a special consideration.   To illustrate these
statements of the problem let us consider the following  process 
\begin{equation}
\label{2}
{\rm d}X_t=-\rho _1X_t\,\1_{\left\{X_t<\vartheta \right\}}{\rm
d}t-\rho _2X_t\,\1_{\left\{X_t\geq \vartheta \right\}}{\rm
d}t+\sigma {\rm d}W_t,\;\; 0\leq t\leq T ,
\end{equation}
where $W_t$ is Wiener process, $\rho _1\not =\rho _2$ and $\rho_i>0$.  We call
it {\it Threshold Ornstein-Uhlenbeck} (TOU) process because it can be
considered as a mixture of two different Ornstein-Uhlenbeck processes with
switching.  If we suppose that $\sigma,\rho _1,\rho _2 $ are
known and $\vartheta \in \Theta =\left(\alpha ,\beta \right)$ is unknown
parameter then we obtain the problem of parameter (threshold) $\vartheta $
estimation. 

It is in some sense similar to TAR \eqref{1} and the link between them can be
clarified by the following consideration. Let us 
consider the discrete time approximation of the process \eqref{2} with
$t_{j}=j\delta,j=1,\ldots,n-1 $, where $\delta =T/n$, then we obtain
$$
X_{t_{j+1}}=\left(1-\rho _1\delta \right)
X_{t_{j}}\;\1_{\left\{X_{t_{j}}<\vartheta  \right\}}+\left(1-\rho
_2\delta \right) X_{t_{j}}\;\1_{\left\{X_{t_{j}}\geq \vartheta
\right\}}+\sigma \;\left[W_{t _{j+1}}-W_{t _{j}}\right]. 
$$
This process coincides with \eqref{1} if we put $X_j=X_{t_{j}}, \varrho
_i=\left(1-\rho_i\, \delta  \right)$ and $\varepsilon _{j+1}=\sigma \;\left[W_{t
_{j+1}}-W_{t _{j}}\right]\;\sim {\cal N}\left(0,\sigma ^2\delta \right),$
i.e., $s^2=\sigma ^2\delta $. Hence, the regression model \eqref{1} is a
discrete time approximation of the TOU process \eqref{2}.

The threshold estimation problems for both models are of singular type and the
limit distributions of the MLE's $n(\hat\vartheta _n-\vartheta
)$ and $T(\hat\vartheta _T-\vartheta
)$ are of $arg sup$ type functionals of the compound Poisson and Wiener
processes respectively.

The process $\left(X_t\right)_{t\geq 0}$ has ergodic properties, the
invariant density  is a mixture of two Gaussian,  
the Fisher information is equal to infinity and we show that  the maximum
likelihood and 
Bayesian estimators converge to two different limit laws. 

We consider  several other  threshold type models of ergodic
diffusion processes and study  the asymptotic properties of  the ML and Bayesian
estimators. We discuss as well the construction of the goodness of fit tests for
such threshold models.

\section{Threshold Ornstein-Uhlenbeck Process}

\subsection{Threshold estimation}

We start with the  TOU process 
\begin{equation}
\label{3}
{\rm d}X_t=-\rho _1X_t\,\1_{\left\{X_t<\vartheta \right\}}{\rm
d}t-\rho _2X_t\,\1_{\left\{X_t\geq \vartheta \right\}}{\rm
d}t+\sigma {\rm d}W_t,\;\;X_0,\;\; 0\leq t\leq T ,
\end{equation}
 where we  suppose that the following condition is fulfilled.

{\bf Condition ${\cal A}$.} {\it The constants $\rho _1\not =\rho _2$,
$\rho_i>0$ and $\sigma ^2>0$ are known and the parameter $\vartheta \in \Theta
=\left(\alpha ,\beta \right), \alpha >0$ is unknown.  The initial value $X_0$ is
independent on the Wiener process random variable.}

The value $\vartheta =0$ is excluded because in the case $\vartheta =0$ there
is no jump in the trend coefficient and the properties of estimators are quite
different. 

 We consider the problem of estimation of the threshold $\vartheta $ by the
continuous time observations $X^T=\left(X_t,0\leq t\leq T\right)$ and we are
interested by the asymptotic behavior of estimators as $T\rightarrow \infty $.

Note that the conditions ${\cal ES}$ of the existence of solution and ${\cal
RP}$ of the ergodicity are fulfilled (see \cite{Kut04}, Sections 1.1 and 1.2)
and the process $\left(X_t\right)_{t\geq 0}$ has ergodic properties with the
invariant density
$$
f\left(\vartheta,x\right)= p_1\left(x,\vartheta \right)\; e^{-\frac{\rho_1\;\left(x
^2-\vartheta ^2\right)}{\sigma 
^2} }+p_2\left(x,\vartheta
\right)e^{-\frac{\rho_2\;\left(x^2-\vartheta ^2\right) }{\sigma ^2} }.
$$
Here  $ p_1\left(x,\vartheta \right)=G\left(\vartheta \right)^{-1}
 \1_{\left\{x< \vartheta \right\}}, $ and $ p_2\left(x,\vartheta
 \right)=G\left(\vartheta \right)^{-1}\1_{\left\{x\geq \vartheta
 \right\}}$ and $G\left(\vartheta \right)$ is the normalizing constant. To
 simplify the exposition we suppose that the random variable $X_0$ has the
 density function $f\left(\vartheta ,x\right)$, hence the observed process is
 stationary. 

We are interested by the asymptotic behavior of the maximum likelihood and
Bayesian estimators of the parameter $\vartheta $, therefore we need the
likelihood ratio function $L\left(\vartheta ,X^T\right) $. This function  can
be written as  (see \cite{LS-01})
\begin{align*}
&\ln L\left(\vartheta ,X^T\right)=-\frac{\rho_1}{\sigma ^2}
\int_{0}^{T}X_t\1_{\left\{X_t< \vartheta \right\}}\,{\rm
d}X_t-\frac{\rho_2}{\sigma ^2}
\int_{0}^{T}X_t\1_{\left\{X_t\geq \vartheta \right\}}\,{\rm
d}X_t\\
&\qquad -\frac{\rho_1 ^2}{2\sigma
^2}\int_{0}^{T}X_t^2\1_{\left\{X_t<\vartheta \right\}}\,{\rm d}t
-\frac{\rho_2 ^2}{2\sigma 
^2}\int_{0}^{T}X_t^2\1_{\left\{X_t\geq \vartheta \right\}}\,{\rm
d}t  +\ln f\left(\vartheta ,X_0\right).
\end{align*}
The contribution of the term $\ln f\left(\vartheta ,X_0\right) $ is
asymptotically negligeable and we will always omitted it for simplicity of
exposition (see the details in \cite{Kut04}).

The MLE $\hat\vartheta _T$ and BE (for quadratic loss function)
$\tilde\vartheta _T$ are defined as usual by the relations
\begin{equation}
\label{bayes}
 L\left(\hat\vartheta _T,X^T\right)= \sup_{\theta \in \Theta }L\left(\theta
 ,X^T\right) \ \mbox{ and } \  \tilde\vartheta _T=\frac{\int_{\alpha }^{\beta }\theta
 \,p\left(\theta \right)L\left(\theta 
 ,X^T\right){\rm d}\theta  }{\int_{\alpha }^{\beta }p\left(\theta \right)L\left(\theta
 ,X^T\right){\rm d}\theta}.
\end{equation}
To describe theirs properties we need the following notations. Let us
introduce
\begin{itemize}
\item the random process
$$
Z_0\left(u\right)=\exp\left\{W\left(u\right)-\frac{\left|u\right|}{2}\right\}
,\qquad u\in \RR,
$$
where $W\left(\cdot \right)$ is  two-sided Wiener process,
\item two random variables
$\hat u$ and    $\tilde u $ defined by the relations
\begin{equation}
\label{4}
 Z_0\left(\hat u \right)=\sup_{u\in\RR} Z_0\left( u \right),\qquad \quad \tilde u
 =\frac{\int_\RR u \,Z_0\left(u \right)\;{\rm d}u }{\int_\RR Z_0\left(u \right)\;{\rm
 d}u}
\end{equation}
\item the function
$$
\Gamma_\vartheta^2=	\frac{\left(\rho_2-\rho_1\right)
^2\,\vartheta ^2}{G\left(\vartheta \right)\,\sigma ^2
}e^{-\frac{\rho
_1^2\vartheta ^2}{\sigma ^2}}.
$$
\end{itemize}

The properties of estimators are given in the following proposition.

\begin{proposition}
\label{T1} Let the condition ${\cal A}$ be fulfilled, then the MLE
$\hat\vartheta _T$ and the BE $\tilde\vartheta _T$ are  uniformly on compacts
 $\KK\subset \Theta $ consistent: for any $\nu >0$
$$
\sup_{\vartheta \in\KK}\Pb_\vartheta \left\{\left|\hat\vartheta _T-\vartheta
\right|>\nu \right\} \longrightarrow 0,
$$ 
 have two  different limit distributions
$$
T\left(\hat\vartheta _T-\vartheta \right)\Longrightarrow \frac{\hat u}{\Gamma
_\vartheta^2} ,\qquad T\left(\tilde\vartheta _T-\vartheta \right)\Longrightarrow
\frac{\tilde u}{\Gamma _\vartheta^2} ,
$$
 theirs moments converge: for any $p>0$
$$
\Ex_\vartheta \left|T\left(\hat\vartheta _T-\vartheta \right)
\right|^p\longrightarrow \Ex\left| \frac{\hat u}{\Gamma _\vartheta^2} \right|^p,\quad 
\Ex_\vartheta \left|T\left(\tilde\vartheta _T-\vartheta \right)
\right|^p\longrightarrow \Ex \left|\frac{\tilde u}{\Gamma _\vartheta^2} \right|^p
$$
\end{proposition}
For the proof see Section 6.

\bigskip

Note that the same normalization and the same type limits (with different
$\Gamma _\vartheta $)    we have in
the problem of delay $\vartheta $ estimation by the observations of the
following Gaussian process
$$
{\rm d}X_t=-\rho\,X_{t- \vartheta} \;{\rm d}t+\sigma {\rm d}W_t,\quad
0\leq t\leq T 
$$
 See details in \cite{KK} (or in \cite{Kut04}, Section 3.3).

\bigskip

Remind that  the Bayesian estimators are usually  asymptotically efficient in
singular parameter estimation problems
\cite{IH81}.  The following lower bound is valid: for all estimators $\bar\vartheta _T$
\begin{equation*}
\Liminf_{\delta \rightarrow 0}\Liminf_{T \rightarrow\infty
}\sup_{\left|\vartheta -\vartheta_0 \right|<\delta }T^2 \,\Ex_\vartheta
\left(\bar \vartheta _T-\vartheta \right)^2\geq 
\frac{\Ex\tilde u^2}{\Gamma{\left(\vartheta_0\right)}^4} 
\end{equation*}
see \cite{IH81}, Section 1.9 (or \cite{Kut04}, Proposition 2.24). We call an estimator
$\vartheta _T^*$ asymptotically efficient if for all $\vartheta _0\in \Theta $
we have the equality
\begin{equation*}
\lim_{\delta \rightarrow 0}\lim_{T \rightarrow\infty
}\sup_{\left|\vartheta -\vartheta _0\right|<\delta }T^2 \,\Ex_\vartheta
\left( \vartheta _T^*-\vartheta \right)^2=
\frac{\Ex\tilde u^2}{\Gamma{\left(\vartheta_0\right)}^4} .
\end{equation*}
It can be verified that the convergence of the moments of bayesian estimators
is uniform on the compacts in $\Theta $ and that the function $\gamma
\left(\vartheta \right)$ is continuous. From these properties we obtain
immediately the asymptotic efficiency of the bayesian estimators (in the sense
of this lower bound). 

The quantities $\Ex  \hat u^2$ and $\Ex  \tilde
u^2$  were calculated by Terent'ev (1968) and Rubin and Song (1995) respectively 
$$
\Ex \hat u^2=26\quad > \quad \Ex \tilde u^2=16 \zeta \left(3\right)\sim 19,3
$$
where $\zeta \left(\cdot \right)$ is Riemann zeta function. This relation
shows the difference between the limit variances of the MLE and BE. 

\bigskip
\subsection{All parameters unknown}
It is possible to describe the properties of estimators in the case when all
three parameters $(\rho _1,\rho _2,
\vartheta)=(\vartheta_1,\vartheta_2,\vartheta_3)=
\boldsymbol{\vartheta}\in\boldsymbol{\Theta} $ are unknown and we observe
\begin{equation}
\label{5}
{\rm d}X_t=-\vartheta  _1X_t\,\1_{\left\{X_t<\vartheta_3 \right\}}{\rm
d}t-\vartheta  _2X_t\,\1_{\left\{X_t\geq \vartheta_3 \right\}}{\rm
d}t+\sigma {\rm d}W_t,\;\; 0\leq t\leq T .
\end{equation}
We have $\boldsymbol{\Theta}=\left(\alpha _1,\beta _1\right)\times
\left(\alpha _2,\beta _2\right)\times\left(\alpha _3,\beta _3\right)$.  Let us
denote by $\xi $ the random variable with the density $f\left(\boldsymbol{\vartheta}
,x\right)$.

\begin{proposition}
\label{P}
 Suppose that $\beta _1<\alpha _2$ and $\alpha _2>0$, then the MLE
$\hat{\boldsymbol{\vartheta}} _T$, BE $\tilde{\boldsymbol{\vartheta}} _T$ are
consistent, have the following limit distributions 
\begin{align*}
\sqrt{T}\left(\hat\vartheta_{1,T}-\vartheta_1 \right)&\Longrightarrow \zeta
_1\sim {\cal N}\left(0, \frac{\sigma ^2}{\Ex_{\boldsymbol{\vartheta}}\;
\xi^2{\rm\1_{\left\{\xi<\vartheta _3\right\}}}}\right),\\
\sqrt{T}\left(\hat\vartheta_{2,T}-\vartheta_2\right)&\Longrightarrow \zeta
_2\sim{\cal N}\left(0, \frac{\sigma ^2}{\Ex_{\boldsymbol{\vartheta}}\;
\xi^2\1_{\left\{\xi\geq \vartheta _3\right\}}}\right),\\
{T}\left(\hat\vartheta_{3,T}-\vartheta_3\right)&\Longrightarrow \frac{\hat
u}{\Gamma^2 _{\boldsymbol{\vartheta}} },\qquad
{T}\left(\tilde\vartheta_{3,T}-\vartheta_3\right) \Longrightarrow \frac{\tilde
u}{\Gamma^2 _{\boldsymbol{\vartheta}} } .
\end{align*} 
{\it The BE $\tilde \vartheta _{1,T},\tilde \vartheta _{2,T}$ have the same
asymptotic properties as $\hat\vartheta _{1,T},\hat \vartheta _{2,T}$, the
random variables $\zeta _1$ and $\zeta _2$ are independent and are independent
of $\hat u, \tilde u$.}
\end{proposition}

\bigskip

The proof see in Section 6.

The construction of the MLE can be slightly simplified by the following
``separation''. 

 The MLE of the first two components can be written as
$$
\hat\vartheta_{1,T}=-\frac{\int_{0}^{T}X_t\,\1_{\left\{X_t<\hat\vartheta
_{3,T}\right\}}\,{\rm
d}X_t}{\int_{0}^{T}X_t^2\,\1_{\left\{X_t<\hat\vartheta
_{3,T}\right\}}\,{\rm d}t} ,\qquad \hat\vartheta_{2,T}=-\frac{\int_{0}^{T}X_t\,\1_{\left\{X_t\geq \hat\vartheta
_{3,T}\right\}}\,{\rm
d}X_t}{\int_{0}^{T}X_t^2\,\1_{\left\{X_t<\hat\vartheta
_{3,T}\right\}}\,{\rm d}t} 
$$
but to study these expressions can be quite difficult because the estimator
$\hat\vartheta _{3,T}$ depends on the whole trajectory $X^T$ and therefore the
random function $X_t\,\1_{\left\{X_t<\hat\vartheta
_{3,T}\right\}}, \;0\leq t\leq T $ depends of the ``future''. Hence the
stochastic integral needs a special treatment. The problem can be simplified
as follows. Let us estimate the parameter $\vartheta _3$ by the first
$X^{\sqrt{T}}=\left\{X_t,0\leq t\leq \sqrt{T}\right\}$ observation and denote
by $\vartheta _{3,\sqrt{T}}^*$ the 
corresponding consistent  estimator.  We suppose that there exists $b>0$ such
that 
\begin{equation}
\label{con}
\Pb_{\boldsymbol{\vartheta }}\left\{\left|\vartheta _{3,\sqrt{T}}^*-\vartheta
_3  \right|>T^{-b}\right\}\longrightarrow 0 
\end{equation}
as $T\rightarrow \infty $. 
Then we define  the estimators
\begin{equation}
\label{es}
\hat\vartheta_{1,T}=-\frac{\int_{\sqrt{T}}^{T}X_t\,\1_{\left\{X_t<\vartheta _{3,\sqrt{T}}^*\right\}}\,{\rm
d}X_t}{\int_{\sqrt{T}}^{T}X_t^2\,\1_{\left\{X_t<\vartheta _{3,\sqrt{T}}^*\right\}}\,{\rm d}t} ,\quad \hat\vartheta_{2,T}=-\frac{\int_{\sqrt{T}}^{T}X_t\,\1_{\left\{X_t\geq \vartheta
_{3,\sqrt{T}}^*\right\}}\,{\rm
d}X_t}{\int_{\sqrt{T}}^{T}X_t^2\,\1_{\left\{X_t\geq \vartheta
_{3,\sqrt{T}}^*\right\}}\,{\rm d}t}. 
\end{equation}
Now the stochastic integrals are well defined and the consistency and
asymptotic normality of these estimators follow from the usual limit theorems,
i.e., we have
$$
\sqrt{T}\left(\hat\vartheta_{1,T}-\vartheta_{1}\right)=-\sigma \frac{\frac{1}{\sqrt{T}}\int_{\sqrt{T}}^{T}X_t\,\1_{\left\{X_t<\vartheta
_{3,\sqrt{T}}^*\right\}}\,{\rm
d}W_t}{\frac{1}{{T}}\int_{\sqrt{T}}^{T}X_t^2\,\1_{\left\{X_t<\vartheta
_{3,\sqrt{T}}^*\right\}}\,{\rm d}t}
$$
with (law of large numbers)
\begin{equation}
\label{lln}
\frac{1}{{T}}\int_{\sqrt{T}}^{T}X_t^2\,\1_{\left\{X_t<\vartheta
_{3,\sqrt{T}}^*\right\}}\,{\rm d}t\longrightarrow \Ex_\theta \;\xi
^2\1_{\left\{\xi <\vartheta _{3}\right\}}
\end{equation}
and (central limit theorem)
$$
\frac{1}{{\sqrt{T}}}\int_{\sqrt{T}}^{T}X_t\,\1_{\left\{X_t<\vartheta
_{3,\sqrt{T}}^*\right\}}\,{\rm d}W_t\Longrightarrow \zeta \sim {\cal
N}\left(0,\Ex_\theta \;\xi ^2\1_{\left\{\xi <\vartheta
_{3}\right\}}\right).
$$
hence
$$
\sqrt{T}\left(\hat\vartheta_{1,T}-\vartheta_{1}\right)\Longrightarrow {\cal
N}\left(0,\frac{\sigma ^2 }{\Ex_\theta \;\xi ^2\1_{\left\{\xi <\vartheta
_{3}\right\}}}\right).
$$
 Note that the independence of the random variables $\zeta _1$ and
$\zeta _2$ follows from the following property of stochastic integral
$$
\Ex_\theta \left(\int_{0}^{T}X_t\,\1_{\left\{X_t<\vartheta
_{3}\right\}}\,{\rm d}W_t\; \int_{0}^{T}X_t\,\1_{\left\{X_t\geq \vartheta
_{3}\right\}}\,{\rm d}W_t\right)=0.
$$
The possibility to simplify the estimation of $\vartheta _3$ we discuss
at the end of the next section.

\subsection{Misspecification}

 Let us return to the initial problem of threshold estimation and suppose that
 the observed process is
\begin{equation}
\label{6}
{\rm d}X_t=-\rho _1X_t\,\1_{\left\{X_t<\vartheta_0 \right\}}{\rm
d}t-\rho _2X_t\,\1_{\left\{X_t\geq \vartheta_0 \right\}}{\rm
d}t+h\left(X_t\right){\rm d}t+\sigma {\rm d}W_t,
\end{equation}
where $h\left(\cdot \right)$ is some unknown function (contamination) and
$\vartheta _0$ is the true value. We assume that the statistician uses this
model without $h\left(\cdot \right)$ (wrong model) and tries to estimate
$\vartheta $, i.e., he (or she) supposes that the observed process is TOU
\eqref{3} and construct, say, the MLE $\hat\vartheta _T$ as if $h\left(\cdot
\right)\equiv 0$. Then he substitutes the observations \eqref{6} (of course,
containing $h\left(\cdot \right) $. Such situation can be considered as
typical for many applied problems, when there is a difference between the
theoretical model and the real data. Remind that in regular case the MLE and
BE are usually not consistent and converge to the value which minimizes the
Kullback-Leibler distance (see \cite{Kut04}, Section 2.6.1).  The
Kullback-Leibler distance in our problem is (suppose for instant that
$\vartheta _0<\vartheta $)
\begin{align*}
&D_{K-L}\left(\vartheta,\vartheta _0 \right)=\Ex_{\vartheta_0}^* \ln \frac{{\rm
d}\Pb_{\vartheta_0} ^* }{{\rm d}\Pb_{\vartheta}}\left(X^T\right) \\ &\quad
=\frac{T}{2\sigma ^2}\Ex_{\vartheta _0}^* \left[\rho _1\xi
\,\left[\1_{\left\{\xi<\vartheta
\right\}}-\1_{\left\{\xi<\vartheta_0 \right\}}\right]+\rho
_2\xi\,\left[\1_{\left\{\xi\geq \vartheta
\right\}}-\1_{\left\{\xi\geq \vartheta_0
\right\}}\right]+h\left(\xi\right)\right]^2\\
&\quad =\frac{T}{2\sigma ^2}\Ex_{\vartheta _0}^* \left[\left(\rho _1-\rho_2\right)\xi
\,\1_{\left\{\vartheta _0<\xi<\vartheta
\right\}}+h\left(\xi\right)\right]^2\
\end{align*}
where $\Ex_{\vartheta _0}^*$ denotes the expectation w.r.t. the measure
$\Pb_{\vartheta_0} ^*$ which corresponds to the process \eqref{6} (we denote
its density as $f_h\left(\vartheta _0,x\right)$). It can be
shown (see \cite{Kut04}, Section 2.6.1) that 
$$
\hat\vartheta _T\longrightarrow \vartheta _*=\arg\inf_{\vartheta \in \Theta
}D_{K-L}\left(\vartheta,\vartheta _0 \right) .
$$
We are interested by the following question: when $\vartheta _*=\vartheta _0$,
i.e., when the MLE is nevertheless consistent? Surprisingly it is possible
even for not too small functions $h\left(\cdot \right)$. Suppose, for
simplicity, that $\vartheta \in \Theta =\left(\alpha ,\beta \right), \alpha
>0$. 

Let us introduce the function
$$
K\left(\vartheta ,\vartheta _0\right)=\begin{cases}
\vphantom{|_{\big)}}  \Ex_{\vartheta _0}^* \left[\left(\rho _1-\rho_2\right)\xi
\,\1_{\left\{\vartheta _0<\xi<\vartheta
\right\}}+h\left(\xi\right)\right]^2, &\hbox{if }\vartheta\geq \vartheta_0 \\
\vphantom{|^{\big)}}   \Ex_{\vartheta _0}^* \left[\left(\rho _2-\rho_1\right)\xi
\,\1_{\left\{\vartheta<\xi<\vartheta _0
\right\}}+h\left(\xi\right)\right]^2, &\hbox{if }\vartheta \leq \vartheta_0 \\
\end{cases}
$$
and suppose that $\rho _2>\rho _1$. Then for $\vartheta >\vartheta _0$ we have
\begin{align*}
K\left(\vartheta ,\vartheta _0\right)&=\left[\int_{-\infty }^{\vartheta_0
}+\int_{\vartheta }^{\infty }
\right]h\left(x\right)^2f_h\left(\vartheta _0,x\right)\;{\rm d}x\\
&\quad+ \int_{\vartheta_0}^{\vartheta } \left[\left(\rho _1-\rho
_2\right)x+h\left(x\right)\right]^2f_h\left(\vartheta _0,x\right)\;{\rm d}x, 
\end{align*}
and
\begin{align*}
\frac{\partial K\left(\vartheta ,\vartheta _0\right)}{\partial \vartheta
}&=
-h\left(\vartheta \right)^2f_h\left(\vartheta _0,\vartheta \right) +\left[\left(\rho _1-\rho
_2\right)\vartheta +h\left(\vartheta\right)\right]^2f_h\left(\vartheta
_0,\vartheta\right)\\
&= \left[\left(\rho _1-\rho_2\right)^2\vartheta^2+2\left(\rho
_1-\rho_2\right)\vartheta h\left(\vartheta \right)\right] f_h\left(\vartheta
_0,\vartheta \right).
\end{align*} 
Therefore, if 
\begin{align*}
h\left(y\right)<\frac{y}{2}\left(\rho _2-\rho _1\right),\quad {\rm for}\quad
 \alpha<y<\beta, 
\end{align*}
then for $\vartheta >\vartheta _0$
$$
\frac{\partial K\left(\vartheta ,\vartheta _0\right)}{\partial \vartheta}>0
$$
and similarly, if 
\begin{align*}
h\left(y\right)>-\frac{y}{2}\left(\rho _2-\rho _1\right),\quad {\rm for}\quad
 \alpha<y<\beta, 
\end{align*}
then for $\vartheta <\vartheta _0$
$$
\frac{\partial K\left(\vartheta ,\vartheta _0\right)}{\partial \vartheta}<0.
$$
We see that if the function $h\left(\cdot \right)$ satisfies the condition
\begin{equation}
\label{7}
\left|h\left(y\right)\right|<\frac{y}{2}\left(\rho _2-\rho
_1\right),\qquad  \alpha <y<\beta ,
\end{equation}
then $\vartheta_*=\vartheta _0 $ and the MLE $\hat\vartheta _T$ is consistent
even for this ``wrong model'' (see \cite{Kut04}, Section 3.4.5 for another
example). Note, that there is no conditions on $h\left(y\right)$ for $y
\not\in \left[\alpha ,\beta \right]$. 

 \bigskip

Let us return to the problem of the construction of the preliminary consistent
estimator of the parameter $\vartheta _3$ by observations \eqref{5}. Suppose
that $\beta _1-\alpha _1<\alpha _2-\beta _1$ and  $\beta _2-\alpha _2<\alpha
_2-\beta _1$. Let us put 
$$
\hat\vartheta _1=\frac{\alpha _1+\beta _1}{2},\qquad \hat\vartheta
_2=\frac{\alpha _2+\beta _2}{2} 
$$
and consider the problem of estimation $\vartheta _3$ by the ``wrong model''
$$
{\rm d}X_t=-\hat\vartheta  _1X_t\,\1_{\left\{X_t<\vartheta_3 \right\}}{\rm
d}t-\hat\vartheta  _2X_t\,\1_{\left\{X_t\geq \vartheta_3 \right\}}{\rm
d}t+\sigma {\rm d}W_t,\;\; 0\leq t\leq \sqrt{T} 
$$
with ``known'' $\hat\vartheta _1,\hat\vartheta _2$.  This corresponds well to
the model \eqref{6} with
$$
h\left(x\right)=(\hat\vartheta _1-\vartheta
_1)x\1_{\left\{x<\vartheta_3 \right\}} +(\hat\vartheta _2-\vartheta
_2)x\1_{\left\{x\geq \vartheta_3 \right\}}.
$$
We see that the condition \eqref{7} is fulfilled, hence the MLE $\hat\vartheta
_{3,\sqrt{T}}$ is consistent and can be used in the construction of the
estimators \eqref{es}. Note that the  estimator $\hat\vartheta
_{3,\sqrt{T}}$  even has ``singular'' rate of
convergence, but its limit distribution is different of that of the true MLE.

\section{Other Threshold Models.}

Below we consider several other threshold type ergodic diffusion processes and
discuss the properties of parameter  estimators for these models.

\subsection{Simple Threshold model.} Suppose that the observed process is
\begin{equation}
\label{9}
{\rm d}X_t=\rho _1\,\1_{\left\{X_t<\vartheta \right\}}{\rm d}t-\rho
_2\,\1_{\left\{X_t\geq \vartheta \right\}}{\rm d}t+\sigma {\rm d}W_t,\quad
0\leq t\leq T ,
\end{equation}
where $\rho _i>0$ and $\vartheta \in \left(\alpha ,\beta \right)$. Then this
process is ergodic with exponential type invariant density
$$
f\left(\vartheta ,x\right)=\frac{1}{G\left(\vartheta \right)}\exp\left\{{-\frac{2\rho
\left(x,\vartheta \right)\left|x-\vartheta \right|}{\sigma ^2}}\right\},
$$
where $\rho \left(x,\vartheta \right)=\rho _1\1_{\left\{x<\vartheta
 \right\}}+\rho _2\1_{\left\{x\geq \vartheta \right\}}$ and $G\left(\vartheta
 \right)$ is the normalizing constant.

The MLE $\hat\vartheta _T$ and BE $\tilde \vartheta _T$ have the same
properties as in Theorem \ref{T1}
$$
T\left(\hat\vartheta _T-\vartheta \right)\Longrightarrow \frac{\hat u}{\Gamma^2
_\vartheta } ,\qquad T\left(\tilde\vartheta _T-\vartheta
\right)\Longrightarrow \frac{\tilde u}{\Gamma^2 _\vartheta }
$$
 and the corresponding function 
$$
\Gamma_\vartheta^2=	\frac{2\,\rho_2\rho_1\left(\rho_2+\rho_1\right)}{\sigma ^4}  .
$$ 

Note, that the normalized LR converges to the limit process as follows:
$$
Z_T\left(u\right)=\frac{L\left(\vartheta +\frac{u}{T},
X^T\right)}{L\left(\vartheta , X^T\right)}\Longrightarrow \exp\left\{\Gamma
_\vartheta \;W\left(u\right)-\frac{\left|u\right|}{2}\;\Gamma 
_\vartheta^2\right\}.
$$
The proof see in the Section 6.

\subsection{Simple Switching.} Suppose that in the model \eqref{9} we have
$\rho _1=\rho _2=\rho >0$. Then the observed process is
\begin{equation}
\label{ssp}
{\rm d}X_t=-\rho\,\sgn \left(X_t- \vartheta \right){\rm d}t+\sigma {\rm d}W_t,\quad
0\leq t\leq T ,
\end{equation}
where $\vartheta \in \Theta =\left(\alpha ,\beta \right)$. This {\it Simple
Switching Process} was studied in \cite{Kut04}, Section 3.4.1. Remind that it has 
Laplace type invariant density
$$
f\left(\vartheta ,x\right)= \frac{\rho }{\sigma ^{2}}\;e^{-\frac{2\rho}{\sigma
^2} \left|x-\vartheta\right|} .
$$ 
The likelihood ratio formula has the representation
$$
L\left(\vartheta ,X^T\right)=\exp\left\{-\frac{\rho}{\sigma ^2} \int_{0}^{T}\sgn \left(X_t-
\vartheta \right){\rm d} X_t-\frac{\rho ^2T}{2\sigma ^2}\right\}.
$$
Hence, the MLE $\hat\vartheta _T$ is defined by the equation
$$
\int_{0}^{T}\sgn \left(X_t- \hat\vartheta _T \right){\rm d}
X_t=\inf_{\vartheta \in\left(\alpha ,\beta \right)}\int_{0}^{T}\sgn \left(X_t-
\vartheta \right){\rm d} X_t.
$$
Note that the last stochastic integral we find in Tanaka-Meyer representation of
the local time of diffusion process (see \cite{RY})
$$
\Lambda _T\left(\vartheta \right)=\left|X_T-\vartheta
\right|-\left|X_0-\vartheta \right|-\int_{0}^{T}\sgn \left(X_t- 
\vartheta \right){\rm d} X_t
$$
and the maximum likelihood is in some sense asymptotically equivalent to the
{\it maximum local time estimator}. Remind that
$f_T^\circ\left(x\right)=\Lambda _T\left(x \right)/T\sigma ^2$ is the
consistent, asymptotically normal and asymptotically efficient (in
nonparametric statement) estimator of
the invariant density (see \cite{Kut04} for details),  and we have obviously
$$
\sup_{\vartheta \in\Theta }f\left(\vartheta _0,\vartheta
\right)=f\left(\vartheta _0,\vartheta _0 \right) .
$$
.

We have the same asymptotic properties of the MLE and BE as in the Theorem \ref{T1}.

The normalized LR
$$
Z_T\left(u\right)=\frac{L\left(\vartheta +\frac{u}{T},
X^T\right)}{L\left(\vartheta , X^T\right)}\Longrightarrow \exp\left\{\Gamma
_\vartheta \;W\left(u\right)-\frac{\left|u\right|}{2}\;\Gamma 
_\vartheta^2\right\},\qquad \Gamma_\vartheta^2=	\frac{4\,\rho^3}{\sigma ^4}  .
$$
 The proof can be found in \cite{Kut04}, Section 3.4.

 \bigskip

The observation window  $\left(-\infty ,\infty \right)$ can be essentially
reduced. Let us put 
$$
\vartheta _{\sqrt{T}}^\star=\frac{1}{\sqrt{T}}\int_{0}^{\sqrt{T}}X_t\;{\rm
d}t.
$$
Note that $\vartheta _{\sqrt{T}}^\star $ is an estimator of the method of
moments ($\Ex_\vartheta \xi =\vartheta $). It is consistent and
asymptotically normal
$$
T^{1/4} \left(\vartheta
_{\sqrt{T}}^\star-\theta \right) \Longrightarrow {\cal
N}\left(0,d^2\left(\vartheta \right)\right),
$$
see \cite{Kut04}, p. 270, where $d\left(\vartheta ^2\right)$ is calculated.
Introduce the {\it window}
$$
\BB_T= \left[\vartheta _{\sqrt{T}}^\star-T^{-1/8} , \vartheta
_{\sqrt{T}}^\star+T^{-1/8} \right].
$$
The MLE and BE we define with the help of the LR $L\left(\vartheta ,X^T_{\sqrt{T}}\right)$
\begin{align*}
=\exp\left\{  -\frac{\rho }{\sigma
^2} \int_{\sqrt{T}}^{T}\sgn{\left(X_t-\vartheta
\right)}\1_{\left\{X_t\in\BB_T\right\}} {\rm d}X_t-\frac{\rho^2
}{2\sigma^2}\int_{\sqrt{T}}^{T}\1_{\left\{X_t\in\BB_T\right\}} {\rm d}t 
\right\}
\end{align*}
Then these estimators have the same
asymptotic  properties as if the observation window is $\BB_T=\left(-\infty ,\infty \right)$. 

This a bit surprising result is probably typical for singular estimation
problems. The analyse of the proof of the properties of estimators (see
\cite{Kut04}, Section 3.4)  shows that only the values of $X_t$ close
to the true value $\vartheta _0$ have contribution to the limit likelihood
ratio. Hence all other observations are irrelevant and can be deleted by
introducing this window.

\subsection{Multy Threshold O-U Process.} Suppose that the observed process is
\begin{equation}
\label{mtou}
{\rm d}X_t=-\sum_{l=1}^{k+1}\rho _l\,X_t\;\1_{\left\{\vartheta _{l-1}<X_t\leq
\vartheta_{l} \right\}}{\rm d}t+\sigma {\rm d}W_t,\quad
0\leq t\leq T ,
\end{equation}
where $\rho _1>0$, $\rho _{k+1}>0$, $\rho _l\not=\rho _m>0$, $\vartheta
_0=-\infty $, $\vartheta _{k+1}=\infty $ and $\boldsymbol{\vartheta} =\left(\vartheta _{1}
,\ldots,\vartheta _{k} \right)\in \boldsymbol{\Theta} =\Theta _1\times \ldots \times \Theta
_k$, $\Theta _l=\left(\alpha _l,\beta _l\right)$, $\beta _l<\alpha
_{l+1}$. Then this process is ergodic and  the normalized likelihood ratio (${\bf
u}=\left(u_1,\ldots,u_k\right)$) has the following limit
$$
Z_T\left({\bf u}\right)=\frac{L\left(\boldsymbol{\vartheta} +\frac{\boldsymbol{u}}{T},
X^T\right)}{L\left(\boldsymbol{\vartheta} , X^T\right)}\Longrightarrow Z\left({\bf
u}\right)=\prod_{l=1}^k\exp\left\{\Gamma 
_l \;W_l\left(u_l\right)-\frac{\left|u_l\right|}{2}\;\Gamma 
_l^2\right\},
$$
where $W_l\left(\cdot \right)$ are independent two-sided Wiener processes. The
estimators $\hat{\boldsymbol{\vartheta}} _T=\left(\hat\vartheta
_{1,T},\ldots,\hat\vartheta _{k,T} \right)$ and $\tilde
{\boldsymbol{\vartheta}} 
_T =\left(\tilde\vartheta _{1,T},\ldots, \tilde\vartheta _{k,T}\right)$ are consistent, have asymptotically 
independent components, 
$$
  T \left(\hat\vartheta _{l,T}-\vartheta _l\right)\Longrightarrow \frac{\hat
  u_l}{\Gamma _l^2},\qquad T \left(\tilde\vartheta _{l,T}-\vartheta
  _l\right)\Longrightarrow \frac{\tilde u_l}{\Gamma _l^2},
$$ 
i.e.; $\left(\hat u_l, \tilde  u_l\right)$  is independent on $\left(\hat u_m,
\tilde  u_m\right)$ if $l\not=m$ and the moments converge. 

The proof see in the section 6.

\section{General Threshold Model.} 

 Suppose that the
observed diffusion process $X^T=\left\{X_t,0\leq t\leq T\right\}$ satisfies
the equation  
\begin{equation}
\label{ntar}
{\rm d}X_t=\sum_{j=1}^{k+1}S _j\left(X_t\right)\,\1_{\left\{\vartheta
_{j-1}<X_t\leq \vartheta_j \right\}}{\rm d}t+\sigma\left(X_t\right) {\rm d}W_t,\quad X_0,
\end{equation}
where $\vartheta _0=-\infty $, $\vartheta _j\in \Theta _j=\left(\alpha
_j,\beta _j\right), j=1,\ldots,k$, $\vartheta _{k+1}=\infty $, $\beta
_j<\alpha _{j+1}$.  The unknown parameter is $\boldsymbol{\vartheta}=
\left(\vartheta _1,\ldots,\vartheta _k\right)\in \boldsymbol{\Theta }=\Theta
_1\times\ldots\times\Theta _k$. Our goal is to estimate $\boldsymbol{\vartheta}$ and to describe the
asymptotic properties of estimators as $T\rightarrow \infty $. As before, we are interested by the
estimators obtained by the Maximum likelihood and Bayesian methods. 

This model can be called ``Nonlinear Threshold Diffusion Process''. Of course,
all considered above models are nonlinear due to the indicator functions. Here
we use the term ``nonlinear'' because the linear function $\rho x$ in the
trend coefficient $-\rho x\,\1_{\left\{\cdot \right\}}$ is replaced by more
general function $S\left( x\right)$.
\bigskip

${\cal ES}$. {\it The functions $S_j(\cdot )$ 
are locally bounded, the function $\sigma \left(\cdot \right)^2$ is continuous
and positive and for some   $A>0$ the condition}
\begin{equation}
\label{esc}
xS_1\left(x\right)\1_{\left\{x<\alpha _1
\right\}}+xS_{k+1}\left(x\right)\1_{\left\{x\geq \beta _k\right\}}+\sigma
\left(x\right)^2 \leq A\left(1+x^2\right)
\end{equation}
{\it holds}

\bigskip

This condition provides the existence of unique weak solution (see \cite{D}).

We suppose that all measures $\left\{{\bf P}_{\boldsymbol{\vartheta}}^{\left(T\right)},
\boldsymbol{\vartheta}\in \boldsymbol{\Theta} \right\}$ induced by this process in the space
$\left({\cal C}\left(0,T\right), {\cal B}\left(0,T\right)\right)$ are
equivalent to the measure ${\bf P} ^{\left(T\right)} $, which corresponds to
the process
$$
{\rm d}X_t=\sigma \left(X_t\right)\,{\rm d}W_t,\qquad X_0,\quad 0\leq t\leq T
$$
(see   \cite{LS-01}). The likelihood ratio 
$$
L\left(\boldsymbol{\vartheta} ,X^T\right)=\frac{{\rm d}{\bf
P}_{\boldsymbol{\vartheta} } ^{\left(T\right)}}{{\rm d}{\bf P}
^{\left(T\right)}}\left(X^T\right),\quad \qquad \boldsymbol{\vartheta}\in
\boldsymbol{\Theta },
$$
in this problem  is the random function 
\begin{align*}
\ln L\left(\boldsymbol{\vartheta} ,X^T\right)=&\sum_{j=1}^{k+1}\int_{0}^{T}\frac{S
_j\left(X_t\right)}{\sigma \left(X_t\right)^2}\,\1_{\left\{\vartheta
_{j-1}<X_t\leq \vartheta_j \right\}}{\rm d}X_t \\
&\qquad\qquad  -
\sum_{j=1}^{k+1}\int_{0}^{T}\frac{S_j\left(X_t\right)^2}{2\sigma
\left(X_t\right)^2}\1_{\left\{\vartheta _{j-1}<X_t\leq \vartheta_j \right\}}{\rm
d}t .
\end{align*}

The MLE $\widehat{\boldsymbol{\vartheta}}_T $ is defined by the same equation
$$
L\left(\widehat{\boldsymbol{\vartheta}}_T ,X^T\right)=\sup_{\boldsymbol{\theta}
\in\boldsymbol{\Theta}  }L\left(\boldsymbol{\theta} ,X^T\right),
$$
where the function $L\left(\boldsymbol{\vartheta} ,X^T\right)$ is not
differentiable with respect to $\boldsymbol{\vartheta}$.

Note that
$$
\1_{\left\{\vartheta _{j-1}<x\leq \vartheta_j \right\}}=\1_{\left\{x\leq
\vartheta_j \right\}}-\1_{\left\{x\leq \vartheta_{j-1} \right\}}.
$$
Hence
\begin{align*}
\sum_{j=1}^{k+1}S_j\left(x\right)\1_{\left\{\vartheta _{j-1}<x\leq
\vartheta_j \right\}}&= \sum_{j=1}^{k+1}S_j\left(x\right)\1_{\left\{x\leq
\vartheta_j \right\}}-\sum_{j=1}^{k+1}S_j\left(x\right)\1_{\left\{x\leq
\vartheta_{j-1} \right\}}\\ 
&=S_{k+1}\left(x\right)+\sum_{j=1}^{k}\left[S_j\left(x\right)-
S_{j+1}\left(x\right)\right]\1_{\left\{x\leq\vartheta_j \right\}}
\end{align*}
and we can write the likelihood ratio as product of $k+1 $ ``likelihood ratios''
\begin{equation}
\label{plr}
\hat{L}\left(\boldsymbol{\vartheta},X^T\right)=\frac{{\rm d}{\bf
P}_{\boldsymbol{\vartheta} } ^{\left(T\right)}}{{\rm d}{\bf P}_0
^{\left(T\right)}}\left(X^T\right)= L_{k+1}\left(X^T\right)\prod_{j=1}^kL_j\left(\vartheta
_j,X^T\right),  
\end{equation}
where
\begin{align*}
\ln L_{k+1}\left(X^T\right)&=\int_{0}^{T}\frac{S_{k+1}\left(X_t\right)}{\sigma
\left(X_t\right)^2}\;{\rm
d}X_t-\int_{0}^{T}\frac{S_{k+1}\left(X_t\right)^2}{2\sigma
\left(X_t\right)^2}\;{\rm d}t
\end{align*}
and
\begin{align*}
\ln L_j\left(\vartheta
_j,X^T\right)&=\int_{0}^{T}\frac{S_j\left(X_t\right)-S_{j+1}\left(X_t\right)}{\sigma
\left(X_t\right)}\,\1_{\left\{X_t\leq \vartheta_j\right\}}\;{\rm d}X_t\\ &\qquad
\qquad
-\int_{0}^{T}\frac{\left[S_j\left(X_t\right)^2-S_{j+1}\left(X_t\right)^2\right]}{2\sigma
\left(X_t\right)^2} \1_{\left\{X_t\leq \vartheta_j\right\}}\;{\rm d}t.
\end{align*}

 This allows us to reduce the calculation of the MLE
$\widehat{\boldsymbol{\vartheta }}_T$ of  
multidimensional parameter $\boldsymbol{\vartheta }$  to $k$ one-dimensional problems :
$$
\hat\vartheta _{j,T}=\argmax_{\vartheta _j\in \Theta _j} L_j\left(\vartheta
_j,X^T\right),\qquad j=1,\ldots,k,
$$
and to put $\widehat{\boldsymbol{\vartheta }}_T=\left(\hat\vartheta
_{1,T},\ldots,\hat\vartheta _{k,T}\right)$.

To introduce the Bayesian estimator $\tilde{\boldsymbol{\vartheta}} _T$ we
suppose that $\boldsymbol{\vartheta}$ is a random vector with a known
continuous positive density {\it a priori} $p\left(\boldsymbol{\theta}
\right),\boldsymbol{\theta}\in \boldsymbol{\Theta} $ and the loss function
$\ell\left(\boldsymbol{u}\right),\boldsymbol{u}\in \RR^k$ is strictly
convex. The estimator $\widetilde{\boldsymbol{\vartheta}} _T$ is defined as
solution of the following equation
$$
\int_{\boldsymbol{\Theta} }^{ }\Ex_{\boldsymbol{\theta}}
\ell\left(\widetilde{\boldsymbol{\vartheta}} _T -\boldsymbol{\theta}
\right)\,p\left(\boldsymbol{\theta} \right){\rm
d}\boldsymbol{\theta}=\inf_{\boldsymbol{\vartheta} \in \boldsymbol{\Theta} }\int_{\boldsymbol{\Theta}
}^{ }\Ex_{\boldsymbol{\theta}} \ell\left(\boldsymbol{\vartheta} -\boldsymbol{\theta}
\right)\,p\left(\boldsymbol{\theta} \right){\rm d}\boldsymbol{\theta}.
$$
Remind that in the case
$\ell\left(\boldsymbol{u}\right)=\left|\boldsymbol{u}\right|^2$ this estimator
is
$$
\widetilde{\boldsymbol{\vartheta}} _T=\frac{\int_{\boldsymbol{\Theta} }^{
}\boldsymbol{\theta} L\left(\boldsymbol{\theta}
,X^T\right)p\left(\boldsymbol{\theta} \right){\rm
d}\boldsymbol{\theta}}{\int_{\boldsymbol{\Theta} }^{ }
L\left(\boldsymbol{\theta} ,X^T\right)p\left(\boldsymbol{\theta} \right){\rm
d}\boldsymbol{\theta}}.
$$
In this case  we can simplify the calculation of the estimator too. 
Suppose that the density 
$p\left(\boldsymbol{\theta}\right)=p_1\left(\theta _1\right) \cdots
p_k\left(\theta _k\right)$ (the components of $\boldsymbol{\vartheta}$ are independent random
variables). Then, using \eqref{plr}, we can write
$$
\tilde{\vartheta }_{j,T}=\frac{\int_{\Theta _j}^{}\theta _j\,L_j\left(\theta
_j,X^T\right)\, p _j\left(\theta _j\right)\,{\rm d}\theta _j}{\int_{\Theta
_j}^{}L_j\left(\theta _j,X^T\right)\, p _j\left(\theta _j\right)\,{\rm
d}\theta _j},\qquad j=1,\ldots,k,
$$
and then to put $\widetilde{\boldsymbol{\vartheta}} _T=\left(\tilde{\vartheta
}_{1,T} ,\ldots,\tilde{\vartheta }_{k,T}\right)$.

The asymptotic behavior of the diffusion process is defined by the following
condition. 

\bigskip

${\bf {\cal A}}$. {\it  The functions   $S_1\left(x\right),\;S_{k+1}(x)$ and  $\sigma
\left(x\right)$  satisfy the conditions
$$
\left|\sigma \left(x\right)\right|^{-1}\leq B\left(1+\left|x\right|^m\right)
$$
with some $B>0$ and $m>0$ and}
\begin{align*}
\Liminf_{x\rightarrow -\infty
}\frac{S_1\left(x\right)}{\sigma \left(x\right)^2}>0,\qquad \quad
\Limsup_{x\rightarrow \infty 
}\frac{S_{k+1}\left(x\right)}{\sigma \left(x\right)^2}<0.
\end{align*}

\bigskip

By this condition the process $\left(X_t\right)_{t\geq 0}$ has ergodic
properties. Let us denote by $f\left(\boldsymbol{\vartheta} ,x\right)$ the density of its
invariant law and by $\xi $ the random variable with such density function. Note
that by this condition $\xi $ has all polynomial moments \cite{Kut04}.

  The {\it identifiability} condition in this statistical problem is the
  following one  
\begin{equation}
\label{id}
\inf_{y \in\left(\alpha_j ,\beta_j \right)}\left|S _j{\left(y
\right)}-S_{j+1}{\left(y\right)}\right|>0,\qquad j=1,\ldots,k.  
\end{equation}

Let us introduce $\hat{\boldsymbol{u}}_{\boldsymbol{\vartheta}}=\left(\hat u_{1,\boldsymbol{\vartheta}
},\ldots,\hat u_{k,\boldsymbol{\vartheta} } \right)$, where  
$$
\hat u_{j,\boldsymbol{\vartheta} }=\frac{\hat u_{j }}{\gamma
_j\left(\boldsymbol{\vartheta} \right)},\quad\gamma
_j\left(\boldsymbol{\vartheta} \right)^2 =
\frac{\left(S_{j+1}\left(\vartheta_j \right)-S_j\left(\vartheta_j \right)
\right)^2}{\sigma\left(\vartheta_j \right) ^2} \;f\left(\boldsymbol{\vartheta}
,\vartheta_j \right),
$$
 and $\hat u_{1 },\ldots,\hat u_{k } $ are independent random variables
 defined by the equalities
$$
\hat u_{j }={\rm argsup}_{u\in \RR}\left[W_j\left(u\right)-\frac{1}{2}\left|u\right|\right].
$$
 Here $W_j\left(\cdot \right),j=1,\ldots,k$ are independent two-sided Wiener processes. 

Let us  define the random vector
$\tilde{\boldsymbol{u}}_{\boldsymbol{\vartheta }}$ as solution of the
following equation
\begin{align*}
\int_{\RR^k}^{}\ell\left(\tilde{\boldsymbol{u}}_{\boldsymbol{\vartheta }}-\boldsymbol{u
}\right)Z\left(\boldsymbol{u}\right){\rm d}\boldsymbol{u }
=\inf_{\boldsymbol{v}\in \RR^k}\int_{\RR^k}^{}\ell\left(\boldsymbol{v}-\boldsymbol{u
}\right)Z\left(\boldsymbol{u}\right){\rm d}\boldsymbol{u },
\end{align*}
where
\begin{equation}
\label{llr}
Z\left(\boldsymbol{u}\right)=\exp\left\{\sum_{j=1}^{k}\left[\gamma
_j\left(\boldsymbol{\vartheta} \right)W_j\left( u_j \right)-\frac{\left|u_j\right|}{2} \gamma
_j\left(\boldsymbol{\vartheta} \right)^2\right] \right\}.
\end{equation}

\begin{theorem}
\label{T2}  Suppose that  these conditions ${\cal ES, A}$ and \eqref{id} are fulfilled, then the MLE
$\widehat{\boldsymbol{\vartheta}} _T$ and bayessian estimator
$\widetilde{\boldsymbol{\vartheta}} _T$ are consistent, have the following
limit distributions:
$$
T\left(\widehat{\boldsymbol{\vartheta}} _T-\boldsymbol{\vartheta}
\right)\Longrightarrow {\hat{\boldsymbol{u}}_{\boldsymbol{\vartheta }}},\qquad \quad
T\left(\widetilde{\boldsymbol{\vartheta}} _T-\boldsymbol{\vartheta}
\right)\Longrightarrow {\tilde{\boldsymbol{u}}_{\boldsymbol{\vartheta }}}
$$
and the moments converge : for any $p>0$
$$
\lim_{T\rightarrow \infty }T^p\;\Ex_\vartheta \left|\widehat{\boldsymbol{\vartheta}} _T-\vartheta
\right|^p ={\Ex \left|{\hat{\boldsymbol{u}}_{\boldsymbol{\vartheta }}}\right|^p},\qquad \lim_{T\rightarrow \infty }T^p\;\Ex_\vartheta \left|\tilde\vartheta _T-\vartheta
\right|^p ={\Ex \left|{\tilde{\boldsymbol{u}}_{\boldsymbol{\vartheta }}}\right|^p}.
$$
\end{theorem}
The proof is given in the  section 6.

\section{Proofs}
 
First note that the parameter estimation problems for the models of the
observations \eqref{3}, \eqref{9} and \eqref{ssp} are particular cases of the
threshold estimation problem for stochastic process \eqref{ntar}. Therefore,
it is sufficient to prove the Theorem \ref{T2}.

The proof of this theorem is based on the two remarkable theorems by Ibragimov
and Khasminskii (\cite{IH81}, Theorems 1.10.1 and 1.10.2) and some results
obtained before in \cite{Kut04}. Let us remind the main steps of this
approach. Introduce the random function (normalized likelihood ratio)
$$
Z_T\left(\boldsymbol{u}\right)=\frac{L\left(\boldsymbol{\vartheta}+
\frac{\boldsymbol{u}}{T},
X^T\right)}{L\left(\boldsymbol{\vartheta} , X^T\right)},\qquad
\boldsymbol{u}\in \boldsymbol{U}_T= U_{1,T}\times\ldots U_{k,T},
$$
where $ U_{j,T}=\left(T\left(\alpha_j -\vartheta_j \right),T\left(\beta_j
-\vartheta_j \right) \right) $. The properties of estimators follow, roughly
speaking, from the weak convergence of this function to the limit random field
\eqref{llr}:
$$
Z_T\left(\boldsymbol{u}\right)=\frac{L\left(\boldsymbol{\vartheta}
+\frac{\boldsymbol{u}}{T}, X^T\right)}{L\left(\boldsymbol{\vartheta} ,
X^T\right)}\Longrightarrow Z\left(\boldsymbol{u}\right).
$$
Suppose that we have already this convergence and (for simplicity) assume that
 $k=1$. Then for the MLE we have ($\vartheta$ is the true value) :
\begin{eqnarray}
&&\Pb_{\vartheta}\left\{T\left(\hat\vartheta _T-\vartheta\right)<x
\right\}=\nonumber\\ 
&&\quad=\Pb\left\{\sup_{T\left(\theta -\vartheta
\right)<x}L\left(\theta,X^T\right) >\sup_{T\left(\theta -\vartheta
\right)\geq x}L\left(\theta,X^T\right) \right\}\nonumber\\ 
&&\quad=\Pb\left\{\sup_{T\left(\theta -\vartheta
\right)<x}\frac{L\left(\vartheta,X^T\right)}{L\left(\vartheta,X^T\right) }
>\sup_{T\left(\theta -\vartheta\right)\geq
x}\frac{L\left(\vartheta,X^T\right)}{L\left(\vartheta_0,X^T\right)
}\right\}\nonumber\\ 
&&\quad=\Pb\left\{\sup_{u<x}Z_T\left(u\right)
>\sup_{u\geq x}Z_T\left(u\right) \right\}\nonumber\\
&&\qquad \qquad \qquad \qquad \longrightarrow \Pb\left\{\sup_{u<x}Z
\left(u\right) >\sup_{u\geq x}Z\left( u\right) \right\}\label{16}\\ 
&&\quad
=\Pb\left(\frac{\hat u }{\gamma{\left(\vartheta\right)^2} }<x\right),\quad {\rm
i.e.}\quad T\left(\hat\vartheta _T-\vartheta \right)\Longrightarrow
\frac{\hat u }{\gamma{\left(\vartheta\right)^2 }}.
\nonumber
\end{eqnarray}
where we put $\theta =\vartheta+T^{-1}u$.

To describe the behavior of the BE we take  we for simplicity  the square
loss function and use the same  change of  variables $\theta =\vartheta
+u/T\equiv \theta _u$ and , 
\begin{eqnarray*}
&& \tilde\vartheta _T=\frac{\int_{\alpha }^{\beta }\theta p\left(\theta
\right)L\left(\theta ,X^T\right){\rm d}\theta }{\int_{\alpha }^{\beta }
p\left(\theta \right)L\left(\theta ,X^T\right){\rm d}\theta } =\vartheta
+\frac{1}{T} \frac{\int_{U_T }u p\left(\theta _u \right)L\left(\theta _u
,X^T\right)\;{\rm d}u }{\int_{U_T } p\left(\theta _u \right){L\left(\theta _u
,X^T\right)}\;{\rm d}u }\\ 
&&=\vartheta+\frac{1}{T} \frac{\int_{U_T }u
p\left(\theta _u \right)\frac{L\left(\theta _u ,X^T\right)}{L\left(\vartheta
,X^T\right)}\;{\rm d}u }{\int_{U_T } p\left(\theta _u
\right)\frac{L\left(\theta _u ,X^T\right)}{L\left(\vartheta ,X^T\right)
}\;{\rm d}u }=\vartheta+\frac{1}{T} \frac{\int_{U_T }u
p\left(\theta _u \right)Z_T\left( u\right)\;{\rm d}u }{\int_{U_T } p\left(\theta _u
\right) Z_T\left(u\right)\;{\rm d}u }
\end{eqnarray*}
 Then, using  the convergence $p\left(\theta
_u \right)\rightarrow 
p\left(\vartheta\right)$,    we can write
\begin{eqnarray}
&&\Pb_{\vartheta}\left\{{T}\left(\tilde\vartheta _T-\vartheta
\right)<x \right\}=\Pb\left\{ \frac{\int_{U_T }u\; p\left(\theta _u
\right)Z_T\left(u\right)\;{\rm d}u }{\int_{U_T } p\left(\theta _u
\right)Z_T\left(u\right)\;{\rm d}u }<x \right\}\nonumber\\ &&\qquad
\qquad \qquad \longrightarrow \Pb\left\{ \frac{\int_{R}u\; Z\left(u \right)\;{\rm d}u
}{\int_{R} Z\left(u \right)\;{\rm d}u }<x \right\} =\Pb\left(\frac{\tilde u
}{\gamma {\left(\vartheta\right)^2}}<x\right).
\label{17}
\end{eqnarray}
The random variables $\hat u$ and $\tilde u$ are defined in \eqref{4}.

We see that to prove the theorem we need to prove the convergences \eqref{16},
\eqref{17}. These convergences together with the estimates on the large
deviations of estimators will provide the convergence of moments. The
corresponding sufficient conditions are given in the mentioned above theorems
by Ibragimov and Khasminskii. Let us introduce the conditions
\begin{description}
 \item[A.] {\it The finite dimensional distributions of the random function
 $Z_T\left(\cdot \right)$ converge to the finite dimensional distributions of
 the function $Z\left(\cdot \right)$}.
\item[B.] {\it There exist constants $B>0, m>0, b>0$ and $d$ such that for any
$R>0$ and $\left|\boldsymbol{u}\right|\leq R,\left|\boldsymbol{v}\right|\leq R $}
\begin{equation}
\label{b}
\Ex_{\boldsymbol{\vartheta}}
\left|Z_T^{\frac{1}{2m}}\left(\boldsymbol{u}\right)-Z_T^{\frac{1}{2m}}
\left(\boldsymbol{v}\right)\right|^{2m}\leq
B\left(1+R^b\right)\left|\boldsymbol{u}-\boldsymbol{v}\right|^d.
\end{equation}
\item[C.] {\it For any $N>0$, there  exists constant $C_N>0$}, such that
\begin{equation}
\label{c}
\Ex_{\boldsymbol{\vartheta}} Z_T^{\frac{1}{2}}\left(\boldsymbol{u}\right)\leq
\frac{C_N}{\left|\boldsymbol{u}\right|^N} 
\end{equation}
 \end{description}

These conditions are the version of the conditions of Theorems 1.10.1 (with $d>k$) and
1.10.2 \cite{IH81}, which we will verify in this work. 

We start with the condition {\bf A}. Let us consider the case when all
$u_j>0$ and denote $h_j\left(x\right)=S_j\left(x\right)/\sigma \left(x\right)$. Note that
\begin{align*}
&\1_{\left\{\vartheta _{j-1}+\frac{u_{j-1}}{T}<X_t\leq \vartheta
_j+\frac{u_j}{T}\right\}} -\1_{\left\{\vartheta _{j-1}<X_t\leq \vartheta
_j\right\}}\\ &\qquad =\1_{\left\{\vartheta _{j}<X_t\leq \vartheta
_j+\frac{u_j}{T}\right\}}-\1_{\left\{\vartheta _{j-1}<X_t\leq \vartheta
_{j-1}+\frac{u_{j-1}}{T}\right\}}=\1_{\left\{\BB_j\right\}}-\1_{\left\{\BB_{j-1}\right\}}
\end{align*}
in obvious notation.

Then the likelihood ratio $Z_T\left(\boldsymbol{u}\right)$ can be
written as follows
\begin{align*}
\ln Z_T\left(\boldsymbol{u}\right)&=\sum_{j=1}^{k+1}\int_{0}^{T}
h_j\left(X_t\right)\left[\1_{\left\{\BB_j\right\}}-\1_{\left\{\BB_{j-1}\right\}}\right]
{\rm d}W_t\\
&\quad -\frac{1}{2}
\sum_{j=1}^{k}\int_{0}^{T}\left[h_j\left(X_t\right)-h_{j+1}
\left(X_t\right)\right]^2\1_{\left\{\BB_j\right\}}{\rm d}t.
\end{align*}
Using the local time estimator $f_T^\circ\left(x\right)$ of the invariant
density $f\left(\boldsymbol{\vartheta },x\right)$  we write
\begin{align*}
&\int_{0}^{T}\left[h_j\left(X_t\right)-h_{j+1}
\left(X_t\right)\right]^2\1_{\left\{\BB_j\right\}}{\rm d}t\\ &\quad
=T\int_{-\infty }^{\infty }\left[h_j\left(x\right)-h_{j+1}
\left(x\right)\right]^2\1_{\left\{\vartheta _j<x\leq \vartheta
_{j}+\frac{u_j}{T} \right\}} \,f_T^\circ\left(x\right) {\rm d}x\\ &\quad
=T\int_{\vartheta _j}^{\vartheta_{j}+\frac{u_j}{T} }
\left[h_j\left(x\right)-h_{j+1} \left(x\right)\right]^2
\,f_T^\circ\left(x\right) {\rm d}x \\ &\quad =T\int_{\vartheta _j }^{\vartheta
_{j}+\frac{u_j}{T} }\left[h_j\left(x\right)-h_{j+1} \left(x\right)\right]^2
\,f\left(\boldsymbol{\vartheta },x\right) {\rm d}x \\ &\quad\quad
+T\int_{\vartheta _j }^{\vartheta _{j}+\frac{u_j}{T}
}\left[h_j\left(x\right)-h_{j+1} \left(x\right)\right]^2
\,\left[f_T^\circ\left(x\right)-f\left(\boldsymbol{\vartheta },x\right)
\right]{\rm d}x.
\end{align*}
For the random function $\eta
_T\left(x\right)=T\left(f_T^\circ\left(x\right)-f\left(\boldsymbol{\vartheta
},x\right) \right)$ we have the estimate: for any $p>0$ there exist constants
$C_*>0$ and $c _*>0$ such that
\begin{equation}
\label{lte}
\Ex_{\boldsymbol{\vartheta}} \left| \eta _T\left(x\right)\right|^p\leq
C_*\,e^{-c_* \left|x\right|}
\end{equation}
see Proposition 1.11 in \cite{Kut04}. This estimate allows us to prove that the
last integral tends to zero as $T\rightarrow \infty $. We have as well
\begin{align*}
&T\int_{\vartheta _j}^{\vartheta_{j}+\frac{u_j}{T} }
\left[h_j\left(x\right)-h_{j+1} \left(x\right)\right]^2
\,f\left(\boldsymbol{\vartheta },x\right) {\rm d}x\\ &\qquad \quad
\longrightarrow u_j \left[h_j\left(\vartheta _j\right)-h_{j+1} \left(\vartheta
_j\right)\right]^2 \,f\left(\boldsymbol{\vartheta },\vartheta
_j\right)=u_j\gamma_j \left(\boldsymbol{\vartheta }\right)^2 .
\end{align*}
Therefore, 
\begin{align*}
\sum_{j=1}^{k}\int_{0}^{T}\left[h_j\left(X_t\right)-h_{j+1}
\left(X_t\right)\right]^2\1_{\left\{\BB_j\right\}}{\rm d}t\longrightarrow
\sum_{j=1}^{k}u_j\gamma_j \left(\boldsymbol{\vartheta }\right)^2  .
\end{align*}
This convergence by the central limit theorem for stochastic integrals yields
the asymptotic normality of the vector $\boldsymbol{\xi}_T =\left(\xi _{1,T},\ldots,\xi _{k,T}\right) $
$$
 \xi
_{j,T}=\int_{0}^{T}\left[h_j\left(X_t\right)-h_{j+1}
\left(X_t\right)\right]\1_{\left\{\BB_j\right\}}{\rm d}W_t\Longrightarrow
{\cal N}\left(0,u_j\gamma_j \left(\boldsymbol{\vartheta }\right)^2 \right)
$$
with asymptotically independent components, because 
$$
\Ex_{\boldsymbol{\vartheta }}\xi _{j,T}\xi _{l,T}=0,\qquad l\not=j.
$$
Moreover, if we put $\xi _{j,T}=\xi _{j,T}\left(u_j\right)$ and consider the vector
$\boldsymbol{\xi}_{{j,T}}=\left(\xi _{j,T}\left(u_{j,1}\right),\right.$ $\left.\ldots,\xi
_{j,T}\left(u_{j,n}\right) \right)$, where $u_{j,1}, \ldots, u_{j,n}$ is some
collection of values from $U_{j,T}$, then 
\begin{align*}
\Ex_{\boldsymbol{\vartheta }}\xi _{j,T}\left(u_{j,r}\right)\xi
_{j,T}\left(u_{j,q}\right)&= T\int_{\vartheta _j}^{\vartheta
_j+\frac{u_{j,r}\wedge
u_{j,q}}{T}}\left[h_j\left(x\right)-h_{j+1}\left(x\right)\right]^2
f\left(\boldsymbol{\vartheta },x\right) {\rm d}x\\
&\qquad \longrightarrow  \left[u_{j,r}\wedge
u_{j,q}\right] \gamma_j \left(\boldsymbol{\vartheta }\right)^2 .
\end{align*}
Using this  equality and preceding limits we can show  the convergence
$$
\left(\xi _{j,T}\left(u_{j,1}\right),\ldots,\xi
_{j,T}\left(u_{j,n}\right) \right)\Longrightarrow \gamma_j
\left(\boldsymbol{\vartheta }\right) \left(W_j\left(u_{j,1}
\right),\ldots,W_j\left(u_{j,1} \right)\right) .
$$
Therefore the condition {\bf A} is fulfilled.

To verify {\bf B} we do it twice. The first time we check this condition with
$m=1$, which is sufficient for Bayes estimators (multidimensional case) and
then (for MLE) we verify it for the partial likelihoods $Z_{j,T}\left(u\right)$.
Following \cite{Kut04}, Lemma 3.28 we write (we suppose that $v_j<u_j$)
\begin{align}
&\Ex_{\boldsymbol{\vartheta
}}\left|Z_T^{1/2}\left(\boldsymbol{u}\right)-Z_T^{1/2}\left(\boldsymbol{v}\right)
\right|^2\leq\frac{1}{4}\sum_{j=1}^{k}
\Ex_*\int_{0}^{T}\left[h_j\left(X_t\right)-h_{j+1}\left(X_t\right)\right]^2
\1_{\left\{\tilde{\BB}_j\right\}}{\rm d}t\nonumber\\ &\quad
=\frac{1}{4}\sum_{j=1}^{k}T\int_{-\infty }^{\infty
}\left[h_j\left(x\right)-h_{j+1}\left(x\right)\right]^2\1_{\left\{\vartheta
_j+\frac{v_j}{T}<x\leq \vartheta _j+\frac{u_j}{T}\right\}}
f_*\left(x\right){\rm d}x \nonumber\\ &\quad =\frac{1}{4}\sum_{j=1}^{k}T\int_{\vartheta
_j+\frac{vj}{T} }^{\vartheta _j+\frac{u_j}{T}
}\left[h_j\left(x\right)-h_{j+1}\left(x\right)\right]^2 f_*\left(x\right){\rm
d}x \nonumber\\ &\quad\leq C\sum_{j=1}^{k} \left|u_j-v_j\right|\leq
C\left\|\boldsymbol{u}-\boldsymbol{v} \right\|.
\label{B1}
\end{align}
Here $\Ex_*$ and $f_*\left(\cdot \right)$ are expectation and invariant
density which correspond to the stochastic differential equation
\begin{align*}
{\rm d}X_t&=\sum_{j=1}^{k+1}S _j\left(X_t\right)\,\left[\1_{\left\{\vartheta
_{j-1}+\frac{u_{j-1}}{T}<X_t\leq \vartheta_j+\frac{u_{j}}{T} \right\}}+\1_{\left\{\vartheta
_{j-1}+\frac{v_{j-1}}{T}<X_t\leq \vartheta_j+\frac{v_{j}}{T}
\right\}}\right]{\rm d}t\\
&\qquad \qquad +\sigma\left(X_t\right) {\rm d}W_t,\qquad X_0,\quad 0\leq t\leq T
\end{align*} 
(see details in \cite{Kut04}, p. 379). The notation $\tilde{\BB}_j$ is clear
from the second line of \eqref{B1}. 

The condition {\bf B} in the case of the study the MLE we check for the
components $Z_{j,T}\left(u_j\right), u_j\in \UU_{j,T}$ separately as follows. Let
us introduce the stochastic process
$$
V_{j,t}=\left(\frac{Z_{j,t}\left(u_j\right)}{Z_{j,t}\left(v_j\right)}\right)^{1/16},\quad
V_{j,0}=1,\quad 0\leq t\leq T 
$$
and denote
$$
g_j\left(x\right)=\frac{S _j\left(x\right)-S _{j+1}\left(x\right)}{\sigma
\left(x\right)}.
$$
Then the process
\begin{align*}
V_{j,t}=\exp&\left\{\frac{1}{16}\int_{0}^{t}\frac{S _j\left(X_s\right)-S
_{j+1}\left(X_s\right)}{\sigma \left(X_s\right)^2}\,\1_{\left\{\vartheta
_{j}+\frac{v_j}{T}<X_s\leq \vartheta_j+\frac{u_{j}}{T} \right\}}\;{\rm
d}X_s\right.\\ &\qquad \qquad \quad \left. -\frac{1}{32}\int_{0}^{t}\frac{S
_j\left(X_s\right)^2-S _{j+1}\left(X_s\right)^2}{\sigma
\left(X_s\right)^2}\,\1_{\left\{\vartheta _{j}+\frac{v_j}{T}<X_s\leq
\vartheta_j+\frac{u_{j}}{T} \right\}}\;{\rm d}s\right\}
\end{align*}
  by It\^o formula admits the representation (under measure
  $\Pb_{\boldsymbol{\vartheta }}^{\left(T\right)}$)
\begin{align*}
V_{j,T}&=1+\frac{1}{16}\int_{0}^{T}V_{j,t}\;\frac{S _j\left(X_s\right)-S
_{j+1}\left(X_s\right)}{\sigma \left(X_s\right)}\,\1_{\left\{\vartheta
_{j}+\frac{v_j}{T}<X_s\leq \vartheta_j+\frac{u_{j}}{T} \right\}}\;\;{\rm
d}W_t\\ &\qquad -\frac{15}{512}\int_{0}^{T}V_{j,t}\;\left(\frac{S
_j\left(X_s\right)-S _{j+1}\left(X_s\right)}{\sigma
\left(X_s\right)}\right)^2\, \1_{\left\{\vartheta _{j}+\frac{v_j}{T}<X_s\leq
\vartheta_j+\frac{u_{j}}{T} \right\}}\; \;{\rm d}t ,
\end{align*}
Remind that
$$
\1_{\left\{\vartheta _{j}+\frac{v_j}{T}<x\leq
\vartheta_j+\frac{u_{j}}{T}
\right\}}\sum_{l=1}^{k+1}S_l\left(x\right)\1_{\left\{\vartheta _{l-1}<x\leq 
\vartheta_l \right\}}= S_{j+1}\left(x\right)\1_{\left\{\vartheta _{j}+\frac{v_j}{T}<x\leq
\vartheta_j+\frac{u_{j}}{T}
\right\}}.
$$
Therefore we can write
\begin{align*}
&\Ex_{\boldsymbol{\vartheta }
}\left|Z_{j,T}^{1/16}\left(u_j\right)-Z_{j,T}^{1/16}\left(v_j\right) \right|^4
=\Ex_{\boldsymbol{\vartheta }} Z_{j,T}^{1/4}\left(v_j\right) \left|1-V_{j,T}
\right|^4 \\ &\quad \leq \left(\Ex_{\boldsymbol{\vartheta }}
Z_{j,T}^{1/2}\left(v_j\right)\right)^{1/2}\left(\Ex_{\boldsymbol{\vartheta }}
\left|1-V_{j,T} \right|^8\right)^{1/2}\leq\left(\Ex_{\boldsymbol{\vartheta }}
\left|1-V_{j,T} \right|^8\right)^{1/2} 
\end{align*}
because $ \Ex_{\boldsymbol{\vartheta }}
Z_{j,T}^{1/2}\left(v_j\right)\leq 1$.  Further
\begin{align}
&\Ex_{\boldsymbol{\vartheta}} \left|1-V_{j,T} \right|^8\leq
C_1\;\Ex_{\boldsymbol{\vartheta}}
\left(\int_{0}^{T}V_{j,t}\;g_j\left(X_t\right)^2\, \1_{\left\{\vartheta
_{j}+\frac{v_j}{T}<X_s\leq \vartheta_j+\frac{u_{j}}{T} \right\}}\; \;{\rm
d}t\right)^8\nonumber\\ &\qquad \quad +C_2\;\Ex_{\boldsymbol{\vartheta}}
\left(\int_{0}^{T}V_{j,t}\;g_j\left(X_t\right)\, \1_{\left\{\vartheta
_{j}+\frac{v_j}{T}<X_s\leq \vartheta_j+\frac{u_{j}}{T} \right\}}\; \;{\rm
d}W_t\right)^8.
\label{int}
\end{align}
For the last (stochastic)  integral we have the  estimates
\begin{align*}
& \Ex_{\boldsymbol{\vartheta}}
\left(\int_{0}^{T}V_{j,t}\;g_j\left(X_t\right)\, \1_{\left\{\vartheta
_{j}+\frac{v_j}{T}<X_s\leq \vartheta_j+\frac{u_{j}}{T} \right\}}\; \;{\rm
d}W_t\right)^8\\ &\qquad \quad \leq C\Ex_{\boldsymbol{\vartheta}}
\left(\int_{0}^{T}V_{j,t}^2\;g_j\left(X_t\right)^2\, \1_{\left\{\vartheta
_{j}+\frac{v_j}{T}<X_s\leq \vartheta_j+\frac{u_{j}}{T} \right\}}\; \;{\rm
d}t\right)^4\\ &\qquad \quad \leq C\Ex_{\boldsymbol{\vartheta}}\sup_{0\leq
t\leq T}V_{j,t}^8 \left(\int_{0}^{T}g_j\left(X_t\right)^2\,
\1_{\left\{\vartheta _{j}+\frac{v_j}{T}<X_s\leq \vartheta_j+\frac{u_{j}}{T}
\right\}}\; \;{\rm d}t\right)^4\\ &\qquad \quad \leq
C\left(\Ex_{\boldsymbol{\vartheta}}\sup_{0\leq t\leq
T}V_{j,t}^{16}\right)^{1/2} \\ &\qquad \quad\qquad \quad \left(\Ex_{\boldsymbol{\vartheta}}
\left(\int_{0}^{T}g_j\left(X_t\right)^2\, \1_{\left\{\vartheta
_{j}+\frac{v_j}{T}<X_s\leq \vartheta_j+\frac{u_{j}}{T} \right\}}\; \;{\rm
d}t\right)^8\right)^{1/2}.
\end{align*}
Remind that $V_t^{16}$ is martingale and
$\Ex_{\boldsymbol{\vartheta}}V_T^{16}=1$. Using once more the  local time estimator of the
density we write 
\begin{align*}
\int_{0}^{T}g_j\left(X_t\right)^2\, \1_{\left\{\vartheta
_{j}+\frac{v_j}{T}<X_s\leq \vartheta_j+\frac{u_{j}}{T} \right\}}{\rm
d}t=T\int_{\vartheta
_{j}+\frac{v_j}{T}}^{\vartheta_j+\frac{u_{j}}{T} }g_j\left(x\right)^2\, f_T^\circ\left(x\right){\rm
d}x.
\end{align*}
Hence
\begin{align*}
&\Ex_{\boldsymbol{\vartheta}}
\left(\int_{0}^{T}g_j\left(X_t\right)^2\, \1_{\left\{\vartheta
_{j}+\frac{v_j}{T}<X_s\leq \vartheta_j+\frac{u_{j}}{T} \right\}}{\rm
d}t\right)^8\\
&\qquad \quad \leq \left(u_j-v_j\right)^7   T
\int_{\vartheta
_{j}+\frac{v_j}{T}}^{\vartheta_j+\frac{u_{j}}{T} }g_j\left(x\right)^{16}\,
\Ex_{\boldsymbol{\vartheta}}f_T^\circ\left(x\right)^8{\rm 
d}x\leq C\left(u_j-v_j\right)^8.
\end{align*}
The expectation $\Ex_{\boldsymbol{\vartheta}}f_T^\circ\left(x\right)^8 $ due
to the estimate \eqref{lte} is a 
bounded function. 
For the first integral in \eqref{int} the similar calculations yield the
estimate
\begin{align*}
\Ex_{\boldsymbol{\vartheta}}
\left(\int_{0}^{T}V_{j,t}\;g_j\left(X_t\right)^2\, \1_{\left\{\vartheta
_{j}+\frac{v_j}{T}<X_s\leq \vartheta_j+\frac{u_{j}}{T} \right\}}\; \;{\rm
d}t\right)^8\leq C\left(u_j-v_j\right)^8.
\end{align*}
Therefore, for $\left|u_j\right|\leq R, \left|v_j\right|\leq R$
\begin{align}
\Ex_{\boldsymbol{\vartheta}}
\left|Z_{j,T}^{1/16}\left(u_j\right)-Z_{j,T}^{1/16}\left(v_j\right)
\right|^8&\leq 
C\left(u_j-v_j\right)^2+\left(u_j-v_j\right)^4\nonumber\\*& \leq
C\left(1+R^2\right)\left|u_j-v_j\right|^2. 
\label{B2}
\end{align}

\bigskip

To verify condition {\bf C} we follow the proof of the Lemmas 3.29 and 2.11 in
\cite{Kut04}.  By condition \eqref{id} we have 
\begin{align*}
&\Ex_{\boldsymbol{\vartheta}}\sum_{j=1}^{k}\int_{0}^{T}\left[h_j\left(X_t\right)-h_{j+1}
\left(X_t\right)\right]^2
\1_{\left\{\vartheta _j<X_t\leq \vartheta _j+\delta _j\right\}}{\rm d}t\\
&\qquad \quad =
T\sum_{j=1}^{k}\int_{\vartheta _j}^{\vartheta _j+\delta _j}\left[h_j\left(x\right)-h_{j+1}
\left(x\right)\right]^2 f\left(\boldsymbol{\vartheta},x\right) \;{\rm d}x\\
&\qquad \quad =T\sum_{j=1}^{k}\kappa _j\delta
_j\left(1+o\left(1\right)\right)\geq \kappa T\left|\boldsymbol{\delta
} \right|
\end{align*}
with some positive constants $\kappa ,\kappa _j$. Here $\boldsymbol{\delta
}=\left(\delta _1,\ldots,\delta _k\right)$ and  we suppose for simplicity that all
$\delta _j>0$. Hence the inequality \eqref{c} follows from the mentioned above
lemmas. 

\bigskip

The properties of BE follow from the Theorem 1.10.2 in \cite{IH81} because
the conditions {\bf A}, \eqref{B1} and \eqref{c} are sufficient for this
theorem.

 For the MLE we do not apply directly the Theorem 1.10.1 in \cite{IH81}
because it requires in condition {\bf B} that $d>k$. We follow the
modification of this theorem discussed in the proof of the Proposition 2.40 in
\cite{Kut04}. Let us consider the vector of likelihood ratios
$\boldsymbol{Y}\left(\boldsymbol{u}\right)_T
=\left(Z_{1,T}^{1/4}\left(u_1\right),\ldots,Z_{k,T}^{1/4}\left(u_1\right)
\right)$. For the components $Z_{j,T}^{1/4}\left(u_j\right),j=1,\ldots,k$ we
have the joint convergence of its  dimensional distributions to the
distribution of the limit random field
$\boldsymbol{Y}\left(\boldsymbol{u}\right)
=\left(Z_{1}^{1/4}\left(u_1\right),\ldots,Z_{k}^{1/4}\left(u_1\right) \right)$
with independent components and the conditions {\bf B} and {\bf C}. Therefore
we have the tightness of the corresponding vector of measures and for each
component we have the large deviations estimates: for any $L>0$ and $N>0$
there exists $C_N>0$ such that
$$
\Pb_{\boldsymbol{\vartheta
}}^{\left(T\right)}\left\{\sup_{\left|u_j\right|>L}Z_{j,T}^{1/4}\left(u_j\right)\geq
\frac{1}{L^N}\right\}\leq \frac{C_N}{L^N}.
$$
These estimates and the  factorization of the likelihood ratio
\eqref{plr} allows us to finish the proof of the properties of MLE mentioned
in Theorem \ref{T2}. Note that the MLE $\hat\vartheta _{j,T}$ can be written as
$$
\hat\vartheta _{j,T}=\argmax_{\theta_j \in \Theta _j}L_j^{1/4}\left(\theta_j ,X^T\right)
$$
too.

\bigskip

To prove the Proposition \ref{P} we consider the normalized likelihood ratio
(we take $u>0$) 
\begin{align*}
\ln
Z_T&\left(v,w,u\right)=\ln\frac{L\left(\vartheta_1+\frac{v}{\sqrt{T}},\vartheta_2+
\frac{w}{\sqrt{T}},\vartheta_3+\frac{u}{{T}},X^T
\right)}{L\left(\vartheta_1,\vartheta_2,\vartheta_3,X^T \right)}\\
&=-\frac{v}{\sigma \sqrt{T}}
\int_{0}^{T}X_t\1_{\left\{X_t<\vartheta _3\right\}}{\rm d}W_t-\frac{w}{\sigma \sqrt{T}}
\int_{0}^{T}X_t\1_{\left\{X_t\geq \vartheta _3\right\}}{\rm d}W_t   \\
&\quad +\left(\vartheta _2-\vartheta
_1+\frac{w-v}{\sqrt{T}}\right)\frac{1}{\sigma}
\int_{0}^{T}X_t\1_{\left\{\vartheta _3< X_t\leq \vartheta
_3+\frac{u}{\sqrt{T}}\right\}}{\rm d}W_t\\
&\quad-\frac{1}{\sigma
^2}\int_{0}^{T}\left[-\frac{v}{\sqrt{T}}\1_{\left\{X_t<\vartheta
_3\right\}}-\frac{w}{\sqrt{T}}\1_{\left\{X_t\geq \vartheta
_3\right\}}\right.\\
&\quad \left.+\left(\vartheta _2-\vartheta
_1+\frac{w-v}{\sqrt{T}}\right)1_{\left\{\vartheta _3< X_t\leq \vartheta
_3+\frac{u}{\sqrt{T}}\right\}} \right]^2X_t^2{\rm d}t\\ 
&\equiv v\Delta _{1,T}+w\Delta _{2,T}+\left(\frac{\vartheta _2-\vartheta
_1}{\sigma }+\frac{w-v}{\sigma \sqrt{T}}\right)\Delta
_{3,T}\left(u\right)-\frac{1}{2}\;J_T, 
\end{align*}
where the last equality introduce the notation for these integrals. 
For the last integral  we can write
\begin{align}
J_T&=\frac{v^2}{\sigma ^2T}\int_{0}^{T}X_t^2 \1_{\left\{X_t<\vartheta
_3\right\}} {\rm d}t+\frac{w^2}{\sigma ^2T}\int_{0}^{T}X_t^2
\1_{\left\{X_t\geq \vartheta _3\right\}} {\rm d}t\nonumber\\ &\quad
+\frac{\left(\vartheta _2-\vartheta _1\right)^2}{\sigma ^2}\int_{0}^{T}X_t^2
\1_{\left\{\vartheta _3<X_t\leq \vartheta _3+\frac{u}{T}\right\}} {\rm
d}t+o\left(1\right).
\label{jT}
\end{align}
For the first two integrals by the law of large numbers we have
\begin{align}
\label{j1}
&\frac{1}{T}\int_{0}^{T}X_t^2 \;\1_{\left\{X_t<\vartheta
_3\right\}}\; {\rm d}t\longrightarrow \;\Ex_{\boldsymbol{\vartheta }}\xi
^2\1_{\left\{\xi <\vartheta 
_3\right\}},\\
&\frac{1}{T}\int_{0}^{T}X_t^2 \;\1_{\left\{X_t\geq \vartheta
_3\right\}}\; {\rm d}t\longrightarrow \;\Ex_{\boldsymbol{\vartheta }}\xi
^2\1_{\left\{\xi \geq \vartheta 
_3\right\}},\label{j2}
\end{align}
and for the last one using the local time estimator of the density we obtain
\begin{align*}
 &\int_{0}^{T}X_t^2\1_{\left\{\vartheta _3<X_t\leq \vartheta
_3+\frac{u}{T}\right\}} {\rm d}t =T\int_{\vartheta _3}^{\vartheta
_3+\frac{u}{T}}x^2\,f_T^\circ\left(x\right)\,{\rm d}x=T\int_{\vartheta
_3}^{\vartheta _3+\frac{u}{T}}x^2\,f\left(\boldsymbol{\vartheta },
x\right)\,{\rm d}x \\ &\qquad +T\int_{\vartheta _3}^{\vartheta
_3+\frac{u}{T}}x^2\,\left(f_T^\circ\left(x\right)-f\left(\boldsymbol{\vartheta
},x\right)\right)\,{\rm d}x =u\,\vartheta _3^2\,f\left(\boldsymbol{\vartheta
},\vartheta _3\right)+o\left(1\right),
\end{align*}
where in $o\left(1\right)$ we used once more the estimate
\eqref{lte}. Therefore
\begin{align*}
J_T\longrightarrow \frac{v^2}{\sigma ^2}\;\Ex_{\boldsymbol{\vartheta }}\xi
^2\1_{\left\{\xi \leq \vartheta 
_3\right\}}+\frac{w^2}{\sigma ^2}\;\Ex_{\boldsymbol{\vartheta }}\xi
^2\1_{\left\{\xi \geq \vartheta 
_3\right\}}+u\,\frac{\left(\vartheta _2-\vartheta _1\right)^2\vartheta
_3^2}{\sigma ^2}\,f\left(\boldsymbol{\vartheta },\vartheta _3\right). 
\end{align*}
For the stochastic integrals $\Delta _{1,T}$ and $\Delta _{2,T}$ from
\eqref{j1}, \eqref{j2} and by
the central limit theorem we have the convergence
\begin{align}
\label{cptj1}
\Delta _{1,T}&\Longrightarrow \zeta _1\;\sim \;{\cal N}\left(0,{\rm
I}_1\right),\qquad {\rm I}_1=\frac{1}{\sigma 
^2}\;\Ex_{\boldsymbol{\vartheta }}\xi ^2\1_{\left\{\xi \leq \vartheta
_3\right\}}\\
\label{cptj2}
\Delta _{2,T}&\Longrightarrow \zeta _2\;\sim \;{\cal N}\left(0,{\rm I}_2
\right),\qquad {\rm I}_2=\frac{1}{\sigma 
^2}\;\Ex_{\boldsymbol{\vartheta }}\xi ^2\1_{\left\{\xi \geq \vartheta
_3\right\}} , 
\end{align}
where the random variables $\zeta _1$ and $\zeta _2$ are independent. 

Let us consider $\Delta _T=\lambda _1\Delta _{3,T}\left(u_1\right)+\lambda _2\Delta
_{3,T}\left(u_2\right)$. We have 
$$
\Delta _T=\int_{0}^{T}\left[\lambda _1 X_t\1_{\left\{\vartheta _3<X_t\leq
\vartheta _3+\frac{u_1}{T}\right\}} +\lambda _2 X_t\1_{\left\{\vartheta _3<X_t\leq
\vartheta _3+\frac{u_2}{T}\right\}}\right]{\rm d}W_t.
$$
Note that
\begin{align*}
&\int_{0}^{T}\left[\lambda _1 X_t\1_{\left\{\vartheta _3<X_t\leq
\vartheta _3+\frac{u_1}{T}\right\}} +\lambda _2 X_t\1_{\left\{\vartheta _3<X_t\leq
\vartheta _3+\frac{u_2}{T}\right\}}\right]^2{\rm d}t\\
&\quad =\lambda _1^2\int_{0}^{T}X_t^2\1_{\left\{\vartheta _3<X_t\leq
\vartheta _3+\frac{u_1}{T}\right\}}{\rm d}t+\lambda
_2^2\int_{0}^{T}X_t^2\1_{\left\{\vartheta _3<X_t\leq 
\vartheta _3+\frac{u_2}{T}\right\}}{\rm d}t\\
&\qquad +2\lambda _1\lambda _2 \int_{0}^{T}X_t^2\1_{\left\{\vartheta _3<X_t\leq
\vartheta _3+\frac{u_1 \wedge u_2}{T}\right\}}{\rm d}t\\
&\longrightarrow \left[u_1\,\lambda _1^2 +u_2\,\lambda _2^2 +2\lambda
_1\lambda _2 \left(u_1 \wedge u_2\right)\right]\vartheta _3^2\;
f\left({\boldsymbol{\vartheta },\vartheta _3 }\right)\equiv d^2.
\end{align*}
Hence $\Delta _T$ is asymptotically normal $\Delta _T\Rightarrow \Delta $ with
the limit variance $d^2$. Remind that the same variance has the random
variable 
$$
\Delta =\lambda _1\vartheta _3 \sqrt{f\left({\boldsymbol{\vartheta },\vartheta
_3 }\right) }\;W\left(u_1\right)+\lambda _2\vartheta _3
\sqrt{f\left({\boldsymbol{\vartheta },\vartheta _3 }\right)
}\;W\left(u_2\right),
$$ 
where $W\left(\cdot \right)$ is a Wiener process. 
Therefore we have the convergence of the finite dimensional distributions of
$\Delta _{3,T}\left(u\right)$ to the
finite dimensional distributions of the process $\vartheta _3
\sqrt{f\left({\boldsymbol{\vartheta },\vartheta 
_3 }\right) }\;W\left(u\right) $: 
\begin{align}
&\left(\Delta _{3,T}\left(u_1\right), \ldots, \Delta
_{3,T}\left(u_k\right)\right)\nonumber\\
&\qquad \Longrightarrow  \left(\vartheta _3
\sqrt{f\left({\boldsymbol{\vartheta },\vartheta 
_3 }\right) }\;W\left(u_1\right), \ldots, \vartheta _3
\sqrt{f\left({\boldsymbol{\vartheta },\vartheta 
_3 }\right) }\;W\left(u_k\right)\right)
\label{kmd}
\end{align}

This convergence together with \eqref{cptj1} and \eqref{cptj2}  allows to
write the likelihood ratio random field as
\begin{align*}
Z_T\left(v,w,u\right)&=\exp\left\{v\Delta _{1,T} -\frac{v^2}{2}{\rm
I}_1+w\Delta _{2,T}-\frac{w^2}{2}{\rm I}_2\right.\\
&\qquad \left. +\left(\frac{\vartheta _2-\vartheta _1}{\sigma }\right)\Delta
_{3,T}\left(u\right) -\frac{\left|u\right|}{2}\gamma
\left(\boldsymbol{\vartheta } \right)^2+ o\left(1\right)\right\}
\end{align*}
where $\Delta _{1,T}  $ and $\Delta _{2,T}  $ are asymptotically normal, and
$$
\gamma
\left(\boldsymbol{\vartheta } \right)^2=\frac{\left(\vartheta _2-\vartheta
_1\right)^2\vartheta _3^2}{ \sigma ^2}\;f\left({\boldsymbol{\vartheta },\vartheta 
_3 }\right) \equiv \gamma ^2.
$$
Therefore we have the convergence of the finite dimensional distributions of
$Z_T\left(v,w,u\right)$ to that of the random function
$$
Z\left(v,w,u\right)=e^{v\zeta _1-\frac{v^2}{2}{\rm I}_1 }\;e^{w\zeta
_2-\frac{w^2}{2}{\rm I}_2 }\; e^{\gamma
W\left(u\right)-\frac{\left|u\right|}{2}\gamma ^2 },\quad v,w,u\in \RR^3
$$
where $\zeta _1$, $\zeta _2$ and $W\left(\cdot \right)$ are independent. 

To check the condition {\bf B} in the case of Bayesian estimation we following
\eqref{B1}  write ($u_2>u_1>0$)
\begin{align*}
&\Ex_{\boldsymbol{\vartheta}}\left|Z_T^{1/2}\left(v_1,w_1,u_1\right)-
Z_T^{1/2}\left(v_2,w_2,u_2\right)\right|^2 \\
 &\quad \leq\frac{1}{4\sigma ^2}
\Ex_*\int_{0}^{T}\left[\frac{\left(v_1-v_2\right)\1_{\left\{X_t<\vartheta
_3\right\}}}{ \sqrt{T}}+\frac{\left(w_1-w_2\right)\1_{\left\{X_t\geq \vartheta
_3\right\}}}{ \sqrt{T}}\right.\\ 
&\qquad \left.+ \left({\vartheta _1-\vartheta
_2}{ }+\frac{v_2-v_1-w_2+w_1}{ \sqrt{T}}\right)\1_{\left\{\vartheta
_3+\frac{u_1}{T}<X_t<\vartheta _3+\frac{u_2}{T}\right\}}\right]^2X_t^2{\rm
d}t\\
 &\quad \leq C_1 \left(v_1-v_2\ \right)^2+C_2 \left(w_1-w_2\
\right)^2+C_3\left|u_2-u_1\right|. 
\end{align*}

In the case of MLE this estimate is not sufficient  because the condition $d>3$
is not fulfilled. We slightly modify the proof of \eqref{B2}. Let us denote
$\boldsymbol{u}=\left(v,w,u\right)$ and put
$$
V_T=\left(\frac{Z_T\left(v_2,w_2,u_2\right)}{Z_T\left(v_1,w_1,u_1\right)}\right)^{\frac{1}{32}}.
$$
Then
\begin{align*}
&\Ex_{\boldsymbol{\vartheta}}\left|Z_T^{\frac{1}{32}}\left(v_1,w_1,u_1\right)-
Z_T^{\frac{1}{32}}\left(v_2,w_2,u_2\right)\right|^8=\Ex_{\boldsymbol{\vartheta}}
Z_T^{\frac{1}{4}}\left(v_2,w_2,u_2\right) \left|1-V_T\right|^8 \\
&\quad \leq
\left(\Ex_{\boldsymbol{\vartheta}}Z_T^{\frac{1}{2}}\left(v_2,w_2,u_2\right)
\right)^{\frac{1}{2}} \left(\Ex_{\boldsymbol{\vartheta}}
\left|1-V_T\right|^{16}  \right)^{\frac{1}{2}}\leq  \left(\Ex_{\boldsymbol{\vartheta}}
\left|1-V_T\right|^{16}  \right)^{\frac{1}{2}} .
\end{align*}
The process $V_t,0\leq t\leq T$ by It\^o formula admits the representation 
$$
V_T=1-a\int_{0}^{T}V_t\left(\Delta S\left(X_t\right)\right)^2{\rm
d}t+b\int_{0}^{T}V_t\left(\Delta S\left(X_t\right)\right){\rm d}W_t 
$$
with corresponding constants $a>0$ and $b>0$ and $\Delta S\left(X_t\right)\equiv \Delta S_t$ is the
difference of two trend coefficients. Hence
\begin{align*}
\Ex_{\boldsymbol{\vartheta}} \left|1-V_T\right|^{16}&\leq
A\Ex_{\boldsymbol{\vartheta}}\left(\int_{0}^{T}V_t\left(\Delta
S_t\right)^2{\rm  d}t\right)^{16}+B\Ex_{\boldsymbol{\vartheta}}\left(\int_{0}^{T}V_t\left(\Delta
S_t\right){\rm  d}t\right)^{16}\\
&\leq A\Ex_{\boldsymbol{\vartheta}}\left(\int_{0}^{T}V_t\left(\Delta
S_t\right)^2{\rm  d}t\right)^{16}+C\Ex_{\boldsymbol{\vartheta}}\left(\int_{0}^{T}V_t^2\left(\Delta
S_t\right)^2{\rm  d}t\right)^{8}.
\end{align*}
Further
\begin{align*}
&\Ex_{\boldsymbol{\vartheta}}\left(\int_{0}^{T}V_t^2\left(\Delta
S_t\right)^2{\rm d}t\right)^{8}\leq \Ex_{\boldsymbol{\vartheta}} \sup_{0\leq
t\leq T}V_t^{16} \left(\int_{0}^{T}\left(\Delta S_t\right)^2{\rm
d}t\right)^{8}\\
 &\quad \leq \left( \Ex_{\boldsymbol{\vartheta}} \sup_{0\leq
t\leq T}V_t^{24} \right)^{\frac{2}{3}}\left(\Ex_{\boldsymbol{\vartheta}}
\left(\int_{0}^{T}\left(\Delta S_t\right)^2{\rm d}t\right)^{24}\right)^{\frac{1}{3}}\\
 &\quad \leq \left(\Ex_{\boldsymbol{\vartheta}}
\left(\int_{0}^{T}\left(\Delta S_t\right)^2{\rm d}t\right)^{24}\right)^{\frac{1}{3}}
\end{align*}
because $\Ex_{\boldsymbol{\vartheta}} \sup_{0\leq
t\leq T}V_t^{24}\leq 1 $. 
Now with the help of \eqref{jT} we can write
\begin{align*}
&\Ex_{\boldsymbol{\vartheta}}\left(\int_{0}^{T}\left(\Delta
S\left(X_t\right)\right)^2{\rm
d}t\right)^{24}=\Ex_{\boldsymbol{\vartheta}}\left(T\int_{-\infty }^{\infty
}\left(\Delta
S\left(x\right)\right)^2f_T^\circ\left(x\right){\rm
d}x\right)^{24}\\ &\qquad \leq C_1\left(v_2-v_1\right)
^{48}+C_2\left(w_2-w_1\right) ^{48}+C_3\left(u_2-u_1\right) ^{24}.
\end{align*}
After substitution of these estimates we obtain 
\begin{align*}
&\Ex_{\boldsymbol{\vartheta}}\left|Z_T^{\frac{1}{32}}\left(v_1,w_1,u_1\right)-
Z_T^{\frac{1}{32}}\left(v_2,w_2,u_2\right)\right|^8\\
&\qquad \quad \leq
A\left|v_2-v_1\right|^8+B\left|w_2-w_1\right|^8+C\left|u_2-u_1\right|^4.
\end{align*}
 Therefore for the values
 $\left|v_i\right|+\left|w_i\right|+\left|u_i\right|\leq R$ we have
\begin{align}
&\Ex_{\boldsymbol{\vartheta}}\left|Z_T^{\frac{1}{32}}\left(v_1,w_1,u_1\right)-
Z_T^{\frac{1}{32}}\left(v_2,w_2,u_2\right)\right|^8\nonumber\\ &\qquad \quad
\leq
C\left(1+R^4\right)\left(\left|v_2-v_1\right|^4+\left|w_2-w_1\right|^4+
\left|u_2-u_1\right|^4\right)
\label{fin}
\end{align}
Hence the condition {\bf B} is fulfilled  with $m=4$ and $d=4>3$ for the random field
$Y_T\left(v,w,u\right)=Z_T^{\frac{1}{4}}\left(v,w,u\right)$. 

To verify the condition {\bf C} we follow the proof of Lemma 2.11 in
\cite{Kut04}. We write ($u>0$)
\begin{align*}
\Ex_{\boldsymbol{\vartheta}}J_T&=\frac{v^2}{\sigma ^2T}\int_{0}^{T}\Ex_{\boldsymbol{\vartheta}}X_t^2 \1_{\left\{X_t<\vartheta
_3\right\}} {\rm d}t+\frac{w^2}{\sigma ^2T}\int_{0}^{T}\Ex_{\boldsymbol{\vartheta}}X_t^2
\1_{\left\{X_t\geq \vartheta _3\right\}} {\rm d}t\\
 &\quad
+\left(\frac{\vartheta _2-\vartheta _1}{\sigma
}+\frac{v-w}{\sigma \sqrt{T}}\right)^2\int_{0}^{T}\Ex_{\boldsymbol{\vartheta}}X_t^2 
\1_{\left\{\vartheta _3<X_t\leq \vartheta _3+\frac{u}{T}\right\}} {\rm
d}t\\
 &\quad+2\frac{w}{\sqrt{T}}\left(\frac{\vartheta _2-\vartheta _1}{\sigma
}+\frac{v-w}{\sigma \sqrt{T}}\right)\int_{0}^{T}\Ex_{\boldsymbol{\vartheta}}X_t^2 
\1_{\left\{\vartheta _3<X_t\leq \vartheta _3+\frac{u}{T}\right\}} {\rm
d}t.
\end{align*}
Note that
\begin{align*}
0<\kappa\equiv\frac{\alpha _2-\beta _1}{\sigma }<\frac{1}{\sigma
}\left|\vartheta _2+\frac{w}{\sqrt{T}}-\left(\vartheta
_1+\frac{v}{\sqrt{T}}\right)\right|< \frac{\beta _2-\alpha _1}{\sigma }\equiv K
\end{align*}
Hence
\begin{align*}
\Ex_{\boldsymbol{\vartheta}}J_T&\geq \frac{v^2}{\sigma
^2}\Ex_{\boldsymbol{\vartheta}}\xi ^2\1_{\left\{\xi <\vartheta _3\right\}}
+\frac{w^2}{\sigma ^2}\Ex_{\boldsymbol{\vartheta}}\xi ^2\1_{\left\{\xi \geq
\vartheta _3\right\}}  +\kappa ^2 T\int_{\vartheta _3}^{\vartheta _3+\frac{u}{T}}
x^2\,f\left(\boldsymbol{\vartheta},x \right){\rm d}x \\ & \qquad
-2\frac{\left|w\right|}{\sqrt{T}}KT\int_{\vartheta _3}^{\vartheta
_3+\frac{u}{T}} x^2\,f\left(\boldsymbol{\vartheta},x \right){\rm d}x.
\end{align*}
Let us put $\delta =\kappa ^2/4K$, then for $\frac{\left|v\right|}{\sqrt{T}}+\frac{\left|w\right|}{\sqrt{T}}+\frac{\left|u\right|}{{T}}\leq \delta 
$ we have
\begin{align}
\Ex_{\boldsymbol{\vartheta}}J_T&\geq v^2{\rm I}_1+w^2{\rm
I}_2+\left|u\right|\frac{\kappa ^2}{2} \alpha _3^2 \inf_{\alpha _3<x\leq \beta
_3}f\left(\boldsymbol{\vartheta},x \right),
\label{fin0}
\end{align}
and for the vector $\boldsymbol{h}=\left(h_1,h_2,h_3\right) $ with  $h_1=\frac{v}{\sqrt{T}},h_2=\frac{w}{\sqrt{T}}, h_3=
\frac{u}{{T}}$,  and $\left\|\boldsymbol{h}\right\| \geq\delta $ we can write
\begin{align}
\frac{\Ex_{\boldsymbol{\vartheta}}J_T}{T}&=h_1^2\int_{-\infty }^{\vartheta
_3}x^2f\left(\boldsymbol{\vartheta},x\right){\rm d}x +h_2^2\int_{\vartheta
_3+h_3 }^{\infty}x^2f\left(\boldsymbol{\vartheta},x\right){\rm d}x\nonumber\\
&\quad 
+\left(\vartheta _1-\vartheta _2+h_1  \right)^2 \int_{\vartheta
_3 }^{\vartheta
_3+h_3}x^2f\left(\boldsymbol{\vartheta},x\right){\rm d}x\nonumber\\
&\geq h_1^2\int_{-\infty }^{\vartheta
_3}x^2f\left(\boldsymbol{\vartheta},x\right){\rm d}x +h_2^2\int_{\beta _3
}^{\infty}x^2f\left(\boldsymbol{\vartheta},x\right){\rm d}x\nonumber\\ 
&\quad 
+\left(\alpha _2-\beta _1 \right)^2 \int_{\vartheta
_3 }^{\vartheta
_3+h_3}x^2f\left(\boldsymbol{\vartheta},x\right){\rm d}x>\kappa _1>0.
\label{fin2}
\end{align}
Here we used the representation
\begin{align*}
&-\vartheta _1\1_{\left\{x<\vartheta _3\right\}}-\vartheta _2\1_{\left\{x\geq
\vartheta _3\right\}}+\left(\vartheta _1+h_1\right)\1_{\left\{x<\vartheta
_3+h_3\right\}}-\left(\vartheta _2+h_2\right)\1_{\left\{x\geq
\vartheta _3+h_3\right\}} \\
&\qquad \quad \quad =h_1 \1_{\left\{x<\vartheta _3\right\}}+h_2 \1_{\left\{x\geq \vartheta
_3+h_3\right\}}+\left(\vartheta _1-\vartheta _2+h_1 \right)\1_{\left\{\vartheta _3< x\leq \vartheta
_3+h_3\right\}}.
\end{align*}

Now having \eqref{fin0}
and \eqref{fin2} we can follow the proof of Lemma 2.11 in \cite{Kut04} and 
obtain the estimate \eqref{c}. 

The properties of modified (simplified) estimators defined by the equalities
\eqref{es} will be proved if we verify the law of large numbers \eqref{lln}. For
any $\varepsilon >0$ using the consistency \eqref{con} we can write
\begin{align*}
&\Pb_{\boldsymbol{\vartheta}}\left\{\left|\frac{1}{\sigma
^2T}\int_{\sqrt{T}}^{T}X_t^2\1_{\left\{X_t<\hat\vartheta _{3,T}\right\}}{\rm
d}t-{\rm I}_1\right| >\varepsilon \right\}\leq\Pb_{\boldsymbol{\vartheta}}\left\{\left|\hat\vartheta
_{3,T}-\vartheta_3 \right|\geq T^{-b} \right\}\\
&\qquad +
\Pb_{\boldsymbol{\vartheta}}\left\{\left|\frac{1}{\sigma
^2T}\int_{\sqrt{T}}^{T}X_t^2\1_{\left\{X_t<\hat\vartheta _{3,T}\right\}}{\rm
d}t-{\rm I}_1\right| >\varepsilon, \left|\hat\vartheta _{3,T}-\vartheta_3
\right|<T^{-b} \right\}\\
&\quad \leq \Pb_{\boldsymbol{\vartheta}}\left\{\sup_{\left|\theta -\vartheta_3 \right|<T^{-b}}\left|\frac{1}{\sigma
^2T}\int_{\sqrt{T}}^{T}X_t^2\1_{\left\{X_t<\theta\right\}}{\rm
d}t-{\rm I}_1\right| >\varepsilon\right\}+o\left(1\right).
\end{align*}
Further
\begin{align*}
&\sup_{\left|\theta -\vartheta_3 \right|<T^{-b}}\left|\frac{1}{\sigma
^2T}\int_{\sqrt{T}}^{T}X_t^2\1_{\left\{X_t<\theta\right\}}{\rm d}t-{\rm
I}_1\right|\\ &\qquad =\sup_{\left|\theta -\vartheta_3
\right|<T^{-b}}\left|\frac{1}{\sigma
^2T}\int_{0}^{T}X_t^2\1_{\left\{X_t<\theta\right\}}{\rm d}t-{\rm
I}_1\right|+o\left(1\right)\\
 &\qquad =\sup_{\left|\theta -\vartheta_3
\right|<T^{-b}}\left|\int_{-\infty }^{\theta }\frac{x^2}{\sigma
^2}f_T^\circ\left(x\right){\rm d}x-\int_{-\infty }^{\vartheta_3
}\frac{x^2}{\sigma ^2}f\left(\boldsymbol{\vartheta}, x\right){\rm d}x\right|+o\left(1\right)\\
&\qquad\leq \left|\int_{-\infty }^{\vartheta_3 }\frac{x^2}{\sigma
^2}\left[f_T^\circ\left(x\right) -f\left(\boldsymbol{\vartheta}, x\right)
\right]{\rm d}x\right|+\int_{\vartheta_3 }^{\vartheta_3+T^{-b} }\frac{x^2}{\sigma
^2}f_T^\circ\left(x\right){\rm d}x+o\left(1\right).
\end{align*}
To finish the proof we just mention, that 
$$
\int_{-\infty }^{\vartheta_3 }\frac{x^2}{\sigma
^2}\left[f_T^\circ\left(x\right) -f\left(\boldsymbol{\vartheta}, x\right)
\right]{\rm d}x\longrightarrow 0
$$
by the law of large numbers.

\section{Goodness of Fit Testing}

Remind two well known goodness of fit (GoF) tests of classical statistics \cite{LR}.  If we
observe $n$ i.i.d. random variables $\left(X_1,\ldots,X_n\right)=X^n$ with
distribution function $F\left(x\right)$ and the basic hypothesis is simple
$$
{\scr H}_0,\qquad F\left(x\right)\equiv 
F_*\left(x\right),\qquad x\in \RR,
$$
 then the Cram\'er-von Mises (C-vM)  and Kolmogorov-Smirnov  (K-S)
 statistics are
$$
\WW_n^2=n\int_{ }^{ }\left[\hat
F_n\left(x\right)-F_*\left(x\right)\right]^2 {\rm d}F_*\left(x\right),\;
\; \DD_n=\sup_x\sqrt{n}\left|\hat F_n\left(x\right)-F_*\left(x\right) \right|
$$
respectively. Here $\hat F_n\left(x\right) $ is the empirical distribution
function.

For continuous $F_*\left(x\right)$ under hypothesis
${\scr H}_0$  we have the convergence 
\begin{align*}
W _n^2\Longrightarrow \int_{0}^{1}W_0\left(s\right)^2{\rm
d}s,\quad 
\DD_n\Longrightarrow \sup_{0\leq s\leq
1}\left|W_0\left(s\right)\right|,
\end{align*}
where $W_0\left(\cdot \right)$ is Brownian bridge. The limit distributions do
not depend on  the model $F_*\left(\cdot \right) $ (the tests are
asymptotically {\it distribution free}) and this essentially simplifies the
choice of the corresponding thresholds for the tests Cram\'er-von Mises  and
Kolmogorov-Smirnov. Note that the both tests are consistent against any fixed
alternative. 

Our goal is to discuss the possibility of the construction of asymptotically
distribution free tests for the mentioned in this work threshold diffusion
processes. 

 Suppose that the basic hypothesis is simple:
$$
{\scr H}_0\qquad :\qquad { the\; observed\; process\;X^T\; is\; TOU
}\left(\vartheta _0\right) 
$$
i.e., the observations $X^T=\left(X_t,0\leq t\leq T\right)$ come from the equation
$$
{\rm d}X_t=-\rho _1X_t\,\1_{\left\{X_t<\vartheta_0 \right\}}{\rm
d}t-\rho _2X_t\,\1_{\left\{X_t\geq \vartheta_0 \right\}}{\rm
d}t+\sigma {\rm d}W_t,\;\; 0\leq t\leq T, 
$$
with known $\vartheta _0$ and we have to test this hypothesis. We propose
below some tests of C-vM and K-S types of asymptotic size $\alpha $.

Let us denote 
$
g\left(x,\vartheta \right)=-\rho _1x\,\1_{\left\{x<\vartheta_0
\right\}}-\rho _2x\,\1_{\left\{x\geq \vartheta_0 \right\}} 
$
and following \cite{DK} introduce the statistics
$$
\WW _T^2=\frac{1}{\sigma
^2T^2}\int_{0}^{T}\left[X_t-X_0-\int_{0}^{t}g\left(X_s,\vartheta
_0\right){\rm d}s\right]^2{\rm d}t,
$$
and
$$
 \DD _T=\frac{1}{\sigma
\sqrt{T}}\sup_{0\leq t\leq T}\left|X_t-X_0-\int_{0}^{t}g\left(X_s,\vartheta
_0\right){\rm d}s\right|
$$
It is easy to see that under ${\scr H}_0$ (= in distribution)
$$
\WW _T^2=\int_{0}^{1}W\left(s\right)^2{\rm d}s,\qquad \quad
\DD_T=\sup_{0\leq s\leq 1}\left|W\left(s\right)\right|,
$$
where $W\left(\cdot \right)$ is Wiener process. Hence the tests 
$$
 \psi_T\left(X^T\right)=\1_{\left\{\WW_T^{2}>c_\alpha
 \right\}},\qquad \quad
 \phi_T\left(X^T\right)=\1_{\left\{\DD_T>d_\alpha \right\}}
$$ 
are distribution free. Here the thresholds $c_\alpha ,d_\alpha $ are solutions
of the equations
$$
\Pb\left\{\int_{0}^{1}W\left(s\right)^2{\rm d}s>c_\alpha \right\}=\alpha
,\qquad \quad \Pb\left\{\sup_{0\leq s\leq
1}\left|W\left(s\right)\right|>d_\alpha \right\}=\alpha.  
$$

The both tests are consistent against any fixed alternative..

 Suppose now that the basic hypothesis is composite:
$$
{\scr H}_0\qquad :\qquad { the\; observed\; process\;X^T\; is\; TOU
}\left(\vartheta\right),\vartheta \in \Theta 
$$
 Let us introduce the statistics
$$
\WW _T^2=\frac{1}{\sigma
^2T^2}\int_{0}^{T}\left[X_t-X_0-\int_{0}^{t}g\left(X_s,\vartheta_T^*
\right){\rm d}s\right]^2{\rm d}t,
$$
and
$$
 \DD _T=\frac{1}{\sigma
\sqrt{T}}\sup_{0\leq t\leq T}\left|X_t-X_0-\int_{0}^{t}g\left(X_s,\vartheta_T^*
\right){\rm d}s\right|,
$$
where $\vartheta _T^*$ is the maximum likelihood or bayesian estimator. Remind, that
$\vartheta _T^*=\vartheta + \frac{ u_T^*}{T}$. Using this singular rate of
convergence of estimator it can be shown that under ${\scr H}_0$ we have the
same limit distributions of the statistics
$$
\WW _T^2\Longrightarrow \int_{0}^{1}W\left(s\right)^2{\rm d}s,\qquad \quad
\DD_T\Longrightarrow \sup_{0\leq s\leq 1}\left|W\left(s\right)\right|.
$$
 Hence the tests 
$
 \psi_T\left(X^T\right)=\1_{\left\{\WW_T^2>c_\alpha
 \right\}}$ and  $
 \phi_T\left(X^T\right)=\1_{\left\{\DD_T>d_\alpha \right\}}
$ 
are asymptotically distribution free. These tests as well are consistent against any
fixed alternative.

The similar limits we have in the case of general model \eqref{ntar}.  For
example, suppose
that $\boldsymbol{\vartheta} _0$ is known  and denote the trend coefficient in
\eqref{ntar} as $S\left(\boldsymbol{\vartheta} _0,X_t\right)$. Then once more
(under hypothesis) 
$$
\WW_T^2= \frac{1}{T^2}\int_{0}^{T} \left(\int_{0}^{t}\frac{\left[{\rm
d}X_s-S\left(\boldsymbol{\vartheta} _0,X_s\right){\rm d}s\right]}{\sigma
\left(X_s\right)}\right)^2{\rm d}t=\int_{0}^{1}W\left(s\right)^2{\rm d}s. 
$$
If the basic hypothesis is composite then the same statistic with
$\boldsymbol{\vartheta} _0 $ replaced by one of the estimators (MLE or BE) has
this last integral as limit (in distribution). 

It is interesting to study the direct analogs of the classical C-vM and K-S
tests. Let us introduce the empirical distribution function and empirical
density (local time estimator) 
$$
\hat F_T\left(x\right)=\frac{1}{T}\int_{0}^{T}\1_{\left\{X_t<x\right\}}\,{\rm
d}t,\qquad \qquad  f_T^\circ\left(x\right)=\frac{\Lambda
_T\left(x\right)}{T\sigma \left(x\right)^2} .
$$  
Then the corresponding C-vM statistics
\begin{align*}
\WW_T^2&=T\int_{-\infty  }^{\infty }\left[\hat F_T\left(x\right)-F\left(\boldsymbol{\vartheta}
_0,x \right)\right]^2 {\rm d}F\left(\boldsymbol{\vartheta}
_0,x \right), \\
\VV_T^2&=T\int_{-\infty  }^{\infty }\left[
f_T^\circ\left(x\right)-f\left(\boldsymbol{\vartheta} 
_0,x \right)\right]^2 {\rm d}F\left(\boldsymbol{\vartheta}
_0,x \right)
\end{align*}
have limits in distribution but these limits are  not  distribution
free \cite{Kut04}. One way to have
asymptotically distribution free statistic was proposed by Negri and
 Nishiyama   \cite{NeNi}.  Another possibility (discussed in
\cite{Kut08}) is to use the weight functions.   Let us illustrate the second approach 
on the statistic
$$
\VV_T^2\left(\boldsymbol{\vartheta}_0\right)=T\int_{-\infty  }^{\infty }
H\left(\boldsymbol{\vartheta}_0,x\right) \left( \hat F _T
\left(x\right) -F\left(\boldsymbol{\vartheta}
_0,x\right)\right)^2{\rm d}F\left(\boldsymbol{\vartheta}_0,x \right)
$$
with weight function
$$
H\left(\boldsymbol{\vartheta}_0,x\right)=\frac{\Psi'
\left(\boldsymbol{\vartheta}_0,x\right)
}{f\left(\boldsymbol{\vartheta}_0,x\right)
\left[F\left(\boldsymbol{\vartheta}_0,x\right)-1\right]^2}\,\,M\left({\Psi
\left(\boldsymbol{\vartheta}_0,x\right)}\right).
$$
where  $M{\left(\cdot \right)}$ is some function providing the finitness of
this integral and 
\begin{align*}
\Psi
\left(\boldsymbol{\vartheta}_0,x\right)&=\int_{-\infty}^{x}
\frac{F\left(\boldsymbol{\vartheta}_0,y\right)^2}{\sigma
\left(y\right)^2{f_0\left(\boldsymbol{\vartheta}_0,y\right)}}{\rm d}y\\
&\qquad \qquad +F\left(\boldsymbol{\vartheta}_0,x\right)^2\int_{
x}^{\infty} \left(\frac{F\left(\boldsymbol{\vartheta}_0,y\right)-1}{F\left(\boldsymbol{\vartheta}_0,x\right)-1} \right)^2
\frac{{\rm d}y}{\sigma \left(y\right)^2{f\left(\boldsymbol{\vartheta}_0,y\right)}}.
\end{align*}

It is shown that if $M\left(s\right)=e^{-s}$ then
$$
\VV_T^2\left(\boldsymbol{\vartheta}_0\right)\Longrightarrow \int_{0}^{\infty
}W\left(s\right)^2\,e^{-s}\;{\rm d}s, 
$$
where $W\left(\cdot \right)$ is a Wiener process, i.e.; we have asymptotically
distribution free test $\hat\psi
_T=\1_{\left\{\VV_T^2\left(\boldsymbol{\vartheta}_0\right)>r_\alpha \right\}}$
\cite{Kut08}. The threshold $r_\alpha $, of course, is solution of the
following equation
$$
\Pb\left\{\int_{0}^{\infty
}W\left(s\right)^2e^{-s}\;{\rm d}s>r_\alpha \right\}=\alpha .
$$
The similar result can be proved for the large class of functions
 $M\left(\cdot \right)$ satisfying the obvious conditions.
  
  In the case of
composite hypothesis we can replace $\boldsymbol{\vartheta}_0 $ by one of the
estimators, say, to use
$\VV_T^2(\widehat{\boldsymbol{\vartheta}}_T)$ and to have the same (distribution free)
limit of  this statistic.


\begin{thebibliography}{99}
 
\bibitem{Ch93} Chan, K.S. (1993). Consistency and limiting distribution of the
LSE of a TAR, {\it  Ann. Stat.} {\bf 21}, 520--533. 

\bibitem{CK08} Chan, N.H.  and Kutoyants Yu.A. (2008) On parameter estimation
of threshold autoregressive models, submitted.
\bibitem{DK} Dachian, S. and Kutoyants, Yu.A. (2007)  On the goodness-of-fit
tests for some continuous time processes, in {\sl Statistical Models and Methods
for Biomedical and Technical Systems}, F.Vonta {\sl et al.} (Eds), Birkh\"auser,
Boston,  395-413.
%
\bibitem{D} Durret, R. ( 1996)  {\sl Stochastic Calculus: A Practical
Introduction}. Boca Raton: CRC Press.
\bibitem{FY03} Fan, J. and Yao, Q. (2003). {\it Nonlinear Time
Series: Nonparametric and Parametric Methods.} Springer, New York.
\bibitem{Han00}  Hansen, B.E. (2000). Sample splitting and threshold estimation,
{\it Econometrica} {\bf 68}, 575--603.
\bibitem{IH81} Ibragimov, I.A.  and Khasminskii, R.Z. (1981). {\it Statistical
Estimation}. Springer, New York.

\bibitem{RS} Rubin, H. and  Song, K.-S. (1995). { Exact computation of the asymptotic
efficiency of maximum likelihood estimators of a discontinuous signal in a
Gaussian white noise}. {\it  Ann. Stat.}{\bf 23}, 732--739.

\bibitem{Ter} Terent'yev, A.S.  (1968).  { Probability distribution of a
time location of an absolute maximum at the output of a synchronized filter,}
{\it Radioengineering and Electronics}, {\bf 13},  4, 652--657.
\bibitem{To90} Tong, H. (1990).  {\em Non-linear Time Series: A Dynamical
Systems Approach.}  Oxford University Press, Oxford.

\bibitem{KLS} Koul, H.L., Qianb, L. and Surgailis, D. (2003)
Asymptotics of M-estimators in two-phase linear regression models. {\sl
Stochastic. Process. Appl.}, 103, 123-154.
%
\bibitem{KK}  K\"uchler, U., Kutoyants, Yu. A.  (2000)    Delay estimation for
some stationary diffusion-type processes. {\sl Scand. J.  Statist.},
{\bf 27}, 3, 405--414.
%
\bibitem{Kut04} Kutoyants, Yu.A. (2004) {\sl Statistical Inference for Ergodic
Diffusion Processes,} Springer, London.
%
\bibitem{Kut08} Kutoyants, Yu.A. (2008) On the goodness-of-fit testing for
ergodic diffusion processes, to appear in  {\sl Journal of Nonparametric
Statistics.}
%
\bibitem{LR} Lehmann, E.L. and Romano, J.P. (2005) {\sl Testing Statistical
Hypotheses.} (3rd ed.) Springer, N.Y.
%
\bibitem{LS-01} Liptser, R.S. and Shiryayev, A.N. (2001) {\sl Statistics of
Random Processes. I,} (2nd ed.) Springer, N.Y.
\bibitem{NeNi} Negri, I. and Nishiyama, Y. (2009) { Goodness of fit test for
 ergodic diffusion processes}. {\it Annals of the Institute of
 Statistical Mathematics}, 61,4, 919-928.
\bibitem{Q} Quandt, R.E., (1958) The estimation of the parameters of a linear
regression system obeying two separate regimes. {\sl J. Amer. Statist. Assoc.}
53, 873-880.
\bibitem{RY} Revuz, D. and Yor, M. (1991) {\sl Continuous Martingales and
 Brownian Motion.} Springer, N.Y.


\bibitem{S} Shreve, S.E. (2004) {\sl Stochastic Calculus for Finance II:
Continuous-Time Models}. Springer, N.Y. 

\bibitem{NY} Yoshida, N.(2009) Polynomial type large deviation inequality and
its applications. To appear in {\sl Annals of the Institute of
Statistical Mathematics}.
\end{thebibliography}
\end{document}